\documentstyle[12pt]{article}
\def\m@th{\mathsurround=0pt }
\def\eqalign#1{\null\,\vcenter{\openup\jot \m@th
   \ialign{\strut\hfil$\displaystyle{##}$&$
      \displaystyle{{}##}$\hfil \crcr#1\crcr}}\,}
\newtheorem{theorem}{Theorem}[section]

\newtheorem{corollary}[theorem]{Corollary}
\newtheorem{lemma}[theorem]{Lemma}

\newtheorem{definition}[theorem]{Definition}

\def\ad{{\rm ad }\ }

\def\adr{{\rm ad}_r\ }

\font\goth=eufm10

\def\gg{\hbox{\goth g}}

\def\gn{\hbox{\goth n}}

\def\gh{\hbox{\goth h}}

\def\gs{\hbox{\goth s}}

\def\gl{\hbox{\goth l}}

\def\gr{\hbox{\goth r}}

\def\tip{\hbox{\rm tip}}

\begin{document}

\title
{ Invariant Differential Operators for Quantum Symmetric Spaces, II  }
\author { Gail Letzter\thanks{supported by NSA grant no. 
MDA904-03-1-0033.
AMS subject classification 17B37}\\
Mathematics Department\\
Virginia Tech\\
Blacksburg, VA 24061\\
letzter@math.vt.edu}

\maketitle
 
\begin{abstract} The two papers in this series   
analyze  quantum invariant differential operators for
quantum symmetric spaces  in the maximally split case.  In this paper, we complete the proof of a quantum version of 
Harish-Chandra's theorem:  There is a Harish-Chandra map which 
induces an  isomorphism between the ring of quantum invariant differential operators 
and a ring of  Laurent polynomial invariants with respect to the dotted action of 
the restricted Weyl group.  We
find a particularly nice basis for the quantum invariant differential 
operators that provides a new interpretation 
of  difference 
operators associated to Macdonald polynomials. Finally, we set the 
stage for a general quantum counterpart to noncompact zonal spherical functions.
    \end{abstract}

\section* {Introduction} 
 
Harmonic analysis on symmetric spaces studies
 invariant differential operators and their joint eigenspaces 
in connection with
  Lie groups.   The discovery of quantum groups in the   1980's inspired 
the  growing subject of harmonic analysis on quantum symmetric 
spaces. In particular, zonal spherical functions on most quantum symmetric spaces 
have  been identified with Macdonald or Macdonald-Koornwinder 
polynomials ([N], [NS], [S], [L3], [DN], [NDS], and  [DS]).  In this paper, the focus is  on   quantum invariant differential 
operators.  We complete the proof begun in [L4] of the quantum version 
of Harish-Chandra's fundamental result identifying invariant 
differential operators with restricted Weyl group invariants.  
Moreover, we exhibit a nice basis for the quantum invariant 
differential operators which correspond to a special commuting family of 
difference operators
associated to Macdonald  polynomials.

Let $\gg$ be  a semisimple Lie algebra and  
let $\check U$ denote the simply connected quantized 
enveloping algebra of $\gg$ over the algebraic closure ${\cal C}$ of 
${\bf C}(q)$ ([Jo, Section 3.2.10]).  Quantum symmetric pairs are defined using a left 
coideal subalgebra $B$ which can be viewed as a quantum analog of 
the enveloping algebra for the subalgebra of $\gg$ fixed by an 
involution $\theta$ ([L1, Section 7]). Our investigation of quantum invariant differential 
operators focuses on the ring $\check U^B$ of $B$ invariant elements 
of $\check 
U$  with respect to the right adjoint action.  (Ultimately, it is 
more accurate to identify the ring of quantum invariant differential 
operators with a quotient of the ring $\check U^B$.)

Let $\Sigma$ denote the restricted 
root system associated to $\gg$ and $\theta$.
The Cartan subalgebra of $\check U$ is a group algebra of a 
multiplicative group $\check T$ isomorphic to the weight 
lattice of the root system of $\gg$.  Let ${\cal A}$ be the 
    subgroup 
of $\check T$ which is the image of   the weight lattice of 
$2\Sigma$ under this isomorphism.       The quantum Harish-Chandra projection map
 ${\cal P}_B$
 is a projection of  $\check U$ onto a slight enlargement of ${\cal 
 C}[{\cal A}]$  defined using a 
quantum Iwasawa decomposition ([L4]). It can be viewed as the function 
which computes the eigenvalues of zonal spherical functions with 
respect to the action of elements in $\check U^B$. 
This paper is the second   of two   which analyze the 
image of $\check U^B$ under the Harish-Chandra map ${\cal P}_B$.   

The 
action of the restricted Weyl group $W_{\Theta}$ on $\Sigma$ induces a dotted action of 
$W_{\Theta}$ on the group algebra ${\cal C}[{\cal A}]$.
Building on work in [L4], we complete the proof of a quantum version of 
Harish-Chandra's theorem [HC, Section 4].

 \medskip
 \noindent
 {\bf Theorem A:}  {\it (Corollary 4.2) The Harish-Chandra map ${\cal P}_B$ induces a surjection
 from $\check U^B$ onto the subring of invariants of ${\cal C}[{\cal 
 A}] $ with respect to the dotted action of $W_{\Theta}$. }
 
 \medskip
 Let $F_r(\check 
 U)$ denote the locally finite part of $\check U$ with respect to the 
 right adjoint action. We construct a linear map ${\cal L}$ from 
 $F_r(\check U)$ to $\check U^B$ which is morally 
 the inverse of the Harish-Chandra map ${\cal P}_B$.  This map can be 
 viewed as a generalization of the    map used to analyze the center in 
 [JL] (see also [Jo, Chapter 7].  Write $\tau$ for the 
   isomorphism between ${\cal A}$ and the weight lattice   of $2\Sigma$
   and let $P^+(2\Sigma)$ denote the semigroup of dominant integral restricted 
   weights.  We show that Theorem A is a consequence of fine
   information concerning ${\cal P}_B(\check U^B)$ obtained in the 
   following proposition.

  \medskip
  \noindent
  {\bf Theorem B:} {\it (Theorem 4.1) For each $\lambda\in P^+(2\Sigma)$, there exists   
  $t_{\lambda}$ in the intersection of $\tau(\lambda)+{\rm Ker}\ {\cal 
  P}_B$ and $F_r(\check U)$ such that
  the top degree term of ${\cal P}_B({\cal L}(t_{\lambda}))$ is a 
       nonzero scalar multiple of $\tau( \lambda)$.  Moreover,
        $\{{\cal P}_B({\cal L}(t_{\lambda}))|\lambda\in 
       P^+(2\Sigma)\}$ forms a basis for ${\cal P}_B(\check U^B)$.}
       
       \medskip

Recall that a symmetric pair of Lie algebras can be written as a 
direct sum of irreducible symmetric pairs.
 Theorem A is established in [L4] for all symmetric pairs except 
those which contain a component of type EIII, EIV, EVII, or EIX.    Thus 
the original motivation for this paper was to show that
  Theorem A holds for the remaining four exceptional types.  However, 
Theorem B, which appears quite technical, has significance which goes 
beyond the particulars of this paper.  We will ultimately   use 
Theorem B to develop the theory  of noncompact 
  zonal spherical functions. Recall that compact quantum zonal 
  spherical functions are  elements of the quantized function algebra 
  of the compact Lie group $G$ associated to $\gg$. Unfortunately, 
  there is not yet a  good definition of the 
 quantized function algebra associated to a  noncompact semisimple Lie 
 group.    As a result, the only noncompact zonal 
  spherical functions that have been analyzed so far  are those on the 
  simplest  noncompact   symmetric 
  space associated to $\gs\gl 2$ ([KS]).   
In a future paper,  Theorem B will be the foundation of    a general algebraic definition of zonal spherical functions;
 an outline follows.    Theorem B implies that 
 that $$\check U^B\oplus{\rm Ker}{\cal P}_B= {\cal C}\tau( 
P^+(2\Sigma)) \oplus {\rm Ker}{\cal P}_B.$$  One   shows that $F_r(\check 
U)$ is a subset  of this direct sum of vector spaces.   
  Given an algebra homomorphism $\Lambda$ from ${\cal C}[{\cal A}]$ 
 to ${\cal C}$,  define a 
linear function
$g_{\Lambda}$ on  $\check U^B$ by  
$g_{\Lambda}(a)=\Lambda({\cal P}_B(a))$ 
for   $a\in 
\check U$. Then  $g_{\Lambda}$ extends to a function on the   direct sum 
(and hence on $F_r(\check U)$) where $g_{\lambda}({\rm Ker}\ {\cal 
P}_B)=0$.
Better yet,    $g_{\Lambda}$ extends to  
a  function on $\check U$ by virture of the basic local finiteness 
theorem ([JL,Theorem 6.4])
and an assumption that its restriction  to $\check U^0$ is 
(dotted) invariant with respect to  $W_{\Theta}$. When $\mu$ in $P^+(2\Sigma)$, 
it 
turns out that $g_{q^{\mu}}$ is  a compact quantum zonal 
spherical function at $\mu$ associated to $B$.  Noncompact quantum zonal spherical functions correspond to 
choices of   $\Lambda$  that are not of the form 
$q^{\mu}$ for $\mu\in P^+(2\Sigma)$.  

The reader may wish to focus instead on completing    the 
identification of compact quantum zonal spherical functions with 
orthogonal polynomials Indeed, compact quantum zonal 
spherical functions associated to standard quantum symmetric pairs 
with reduced restricted root systems are identified with Macdonald 
polynomials in [L3] by  computing   radial components 
components of ``small'' elements in $\check U^B$. For most 
irreducible symmetric pairs, these small elements are actually 
contained in the center $Z(\check U)$.   However, this is not true 
for the four problematic exceptional types EIII, EIV, EVII, and EIX.   
In [L3], an  elementwise 
computation was used to find the necessary elements of $\check U^B$ 
for the last three of these types.  
The arguments of this paper provide a simpler proof of their 
existence.  Moreover,  Theorem A 
guarantees the   existence of  $B$ invariant elements 
necessary to make the methods of [L3] extend to the remaining 
problematic type, 
EIII.

Let ${\cal X}$ be the radial 
	   component function on $\check U^B$.   As a  consequence of 
	   Theorem B and properties of ${\cal L}$,
	   the set $\{{\cal X}({\cal L}(t_{\lambda}))|\lambda\in 
	   P^+(2\Sigma)\}$ is a special commuting family of difference 
	   operators associated to Macdonald and related orthogonal 
	   polynomials. (See the discussion following Corollary 4.2 and  
	   [K, Theorem 6.6]).  Thus, Theorem B  leads to a natural 
	   quantum setting for these operators.

We briefly describe  the contents of the rest of this paper. Section 1 sets 
notation  and includes the definitions of the Harish-Chandra map 
${\cal P}_B$ as well as the map ${\cal X}$ which computes   radial 
components.  Results from [L4] are recalled in 
Theorem 1.3.
In Section 2, we establish a 
criterion for ${\cal P}_B(\check U^B)$ to be invariant under the 
dotted
action of $W_{\Theta}$. This involves a careful filtration analysis of ${\cal 
P}_B(Z(\check U))$ and its relation to ${\cal P}_B(\check U^B)$.
 
The map ${\cal L}$ is introduced and studied in Section 3. A formula 
   relating the value of a zonal spherical function at an element 
 $a$ in $F_r(\check U)$ and 
 ${\cal P}_B({\cal L}(a))$ is presented in Lemma 3.5. A study 
   of 
 the possible highest degree term of
 ${\cal L}(a)$ for any $a$ in the set   $ F_r(\check U)\cap (\tau(\lambda)+{\rm 
 Ker}\ {\cal P}_B)$ (Lemma 3.7) shows that
  Theorem A is a consequence of Theorem B (Theorem 3.8).

Theorem B is proved in Section 4.   Using results of earlier 
sections, the proof is reduced to finding a set $\{t_{\lambda}|\lambda\in 
P^+(2\Sigma)\}$ such that each $t_{\lambda}$ lies in the intersection of $\tau(\lambda)+{\rm Ker}\ {\cal 
  P}_B$ and $F_r(\check U)$.  
This is immediate for some
symmetric (see Section 4, Case 
(i)). For   other
symmetric pairs, the proof has two flavors.  We
  use  an induction argument   
exploiting  fine information about the image of central 
elements under ${\cal P}_B$ developed in [L3] for irreducible symmetric 
pairs containing a subpair of type AII (see Section 4, Case (ii)).
The proof for the remaining four types is inspired by the ``small 
element'' arguments 
in [L3, Section 7]. Recall that $F_r(\check U)$ is a direct sum of   
modules  $(\adr U)t$, where $t$ is a group-like element of $\check U$  
and  $(\adr U)t$ contains a unique (up to   scalar) 
nonzero central element.   
Theorem B is established for irreducible symmetric pairs  of type EIII, EIV, EVII, and EIX, 
 by 
showing that the central element of $(\adr 
U)t$ is not a scalar multiple of the $B$ invariant 
element ${\cal L}(t)$ for suitable choices of $t$ (see Section 4,  
Cases (iii) and (iv)).  Rather than dealing with element 
manipulation, we reduce the proof to  
 computations  in the commutative
 polynomial ring ${\cal 
C}[{\cal A}]^{W_{\Theta}.}$ which resemble character formula arguments.

    \section{Background and Notation }
    Let ${\bf C}$ denote the complex numbers, ${\bf Q}$ denote the 
rational numbers,  ${\bf R}$ denote 
the real numbers,  and $q$ an indeterminate. 
Write ${\cal C}$ for the algebraic closure of ${\bf C}(q)$ and let 
${\cal R}$ be the real algebraic closure of ${\bf R}(q)$.

Suppose that $\Phi$ is a root system with $\Phi^+$ equal to the subset 
of positive roots.   Write $Q(\Phi)$ for the root lattice and 
$P(\Phi)$ for the weight lattice associated to $\Phi$.   Set 
$Q^+(\Phi)$ equal to the  span of the positive roots 
$\Phi^+$ over the nonnegative integers.   Let $P^+(\Phi)$ denote the set of dominant integral 
weights associated to $\Phi$.  If a symbol for the set of (positive) simple 
roots associated to $\Phi$ is specified, then we will often replace 
$\Phi$ with this symbol in the notation for the root and weight 
lattices and their subsets. 

Let $\gg$ be a complex semisimple Lie 
algebra.  Let $\gh$ be a fixed Cartan subalgebra of $\gg$ and let 
$\gg=\gn^-\oplus\gh\oplus \gn^+$ be a fixed triangular decomposition.
Let $\Delta$ denote the root system of $\gg$ and 
write $\pi=\{\alpha_1,\dots, \alpha_n\}$  for the  set of positive simple roots.  Let $(\ ,\ )$ denote the Cartan inner product 
on $\gh^*$ associated to $\Delta$.   Recall the standard partial order 
on $\gh^*:  \lambda\geq \gamma$ if and only if 
$\lambda-\gamma\in Q^+(\pi)$.

Let   $\theta$ be a maximally split involution  with respect to the above 
triangular decomposition of $\gg$ 
(see [L1, (7.1), (7.2), and (7.3)]). Write $\gg^{\theta}$ for the corresponding fixed Lie subalgebra of 
$\gg$.
We assume throughout the paper that $\gg,\gg^{\theta}$ is an irreducible 
symmetric pair in the sense of [A] (see also [L2, Section 7]).
It should be noted that the results of this paper extend in a 
straightforward manner to   
general symmetric pairs $\gg,\gg^{\theta}$. 

The involution
  $\theta$ induces an involution $\Theta$ on $\gh^*$ which restricts 
  to an involution on  
  $\Delta$.  Set $\pi_{\Theta}=\{\alpha_i|\ 
  \Theta(\alpha_i)=\alpha_i\}$.
Given $\alpha\in \gh^*$, set $\tilde\alpha= (\alpha-\Theta(\alpha))/2$.  
The restricted root system  $\Sigma$ is the set 
$$\Sigma=\{\tilde\alpha|\alpha\in \Delta{\rm \  and \ }\Theta(\alpha)\neq
\alpha\}.$$  Let $W_{\Theta}$ denote the Weyl group associated to the 
root system $\Sigma$.  

  Let $U=U_q(\gg)$ denote the quantized enveloping algebra of $\gg$ 
which is generated as an algebra   over ${\cal C}$    by $x_i,y_i,t_i^{\pm 1}$ for $1\leq i\leq 
n$.  (See [L1, Section 1, (1.4)-(1.10)] or [Jo, 3.2.9] for the relations and Hopf algebra structure.)  Let 
$T$ denote the group generated by $t_i$ for $1\leq i\leq n$ and let 
$U^0$ denote the subalgebra of $U$ generated by $T$. Let $U^+$
denote the subalgebra of $U$ generated by $x_i, 1\leq i\leq n$, and let
$G^-$ denote the subalgebra of $U$ generated by $y_it_i, 1\leq i\leq n$. 
Recall that there is a group isomorphism
$\tau$  from the additive group 
$Q(\pi)$ to the multiplicative group $T$ defined by 
$\tau(\alpha_i)=t_i$ for $1\leq i\leq n$.

Sometimes it will be necessary to use larger Hopf algebras than $U$ 
which are formed by enlarging $T$.   In particular, suppose that $M$ 
is a
multiplicative monoid  isomorphic  to an  additive submonoid of 
$\sum_{1\leq i\leq n}{\bf Q}\alpha_i$ via the obvious extension of 
$\tau$.   Then $UM$ is the Hopf algebra generated as an algebra by $U$ and $M$ as 
explained in [L4, Section 1]. The most common   extension of $U$ by this method 
is the simply connected quantized enveloping algebra, denoted by 
$\check U$.     As an algebra,  
 $\check U$    is generated 
by $U$ and   $\check T$ where $P(\pi)$ is isomorphic to $\check T$ 
via $\tau$  ([Jo, Section 3.2.10]).  Set $\check U^0$ equal to the 
group algebra of $\check T$ over ${\cal C}$.
Given a subalgebra $A$ of $\check U$, we write $A_+$ for the 
intersection of $A$ with the augmentation ideal of $\check U$.

We associate to $\gg,\gg^{\theta}$ a set   ${\cal B}$ of left coideal 
subalgebras of $U$.  More precisely,  ${\cal B}$   equals   the orbit of  
  real analogs of $U(\gg^{\theta})$ in $U_q(\gg)$ under the 
Hopf algebra automorphisms of $U$ fixing elements of $T$.  The reader 
is referred to [L4, Section 1] for more details. Notation concerning ${\cal B}$ will be 
defined as needed.

Given a multiplicative group $G$, we write ${\cal C}[G]$ for the 
   group algebra generated by $G$ over ${\cal C}$. 
  For any additive group $H$,   write ${\cal C}[H]$ for 
 the group ring generated by $z^{\lambda}$ for $\lambda\in H$. The 
 first notation will be applied to groups  related to $\check T$ while 
 the second 
 notation will be used for groups   related to  $P(\pi)$.
   
Let ${\rm ad}_r$ denote the right adjoint action of $U$ on $\check U$ and 
let ${\rm ad}$ denote the left adjoint action of $U$ on $\check U$.
Given $B$ in ${\cal B}$, write $\check U^B$ for the subset of $B$ invariants in $\check U$ with 
respect to the right adjoint action.  There are two important maps 
     which are used to 
   study $\check U^B$.
The first is the quantum Harish-Chandra projection map and the second 
is the function which computes the quantum radial components.   We 
review the definition of both of these maps here. 

Set $$\check {\cal A}=\{\tau(\tilde\mu)|\mu\in
P(\pi)\}.$$   
Let 
$\check T_{\Theta}=\{\tau((\mu+\Theta(\mu))/2)|\ \mu\in P(\pi)\}$ and 
let ${\cal M}$ denote the subalgebra of $U$ generated by 
$x_i,y_i,t_i^{\pm 1}$ for $\alpha_i\in \pi_{\Theta}$. Let 
$N^+$ be the subalgebra of $U^+$ generated by the elements $(\ad 
x)x_j$ with $\alpha_j\in \pi\setminus\pi_{\Theta}$ and $x\in {\cal 
M}\cap U^+$.  As in 
[L4, (2.4)], we have the following inclusion:
$$\check U\subseteq ((B\check T_{\Theta})_+\check U+N^+_+\check {\cal A})\oplus 
{\cal C}[\check A].\eqno{(1.1)}$$   
\begin{definition} The (quantum) {\it Harish-Chandra map} with respect to 
the symmetric pair $\gg,\gg^{\theta}$ and subalgebra $B$ in ${\cal B}$ is 
the projection   ${\cal P}_B$  of $\check U$ onto ${\cal C}[\check{\cal 
A}]$ using the direct sum  
decomposition (1.1).
\end{definition}

 Set 
 ${\cal A}$ equal to  the group consisting of the elements $\tau(2\tilde\mu)$ 
for $\mu\in P(\Sigma)$.    Write ${\cal C}[Q( \Sigma)]{\cal A}$ for the 
 subring of ${\rm End}_r{\cal C}[P(\Sigma)]$ 
generated by ${\cal C}[Q(\Sigma)]$ and ${\cal A}$ as explained in [L4, 
Section 4].  Given $f\in  {\cal C}[Q( 
\Sigma)]{\cal A}$  and $a\in {\cal C}[Q(\Sigma)]$, we write $a*f$ for 
the right action of $f$ on $a$ ([L4, Section 4]). Let 
 ${\cal C}(Q(\Sigma)){\cal A}$ denote the 
 localization of ${\cal C}[Q(\Sigma)]{\cal A}$ at the Ore set 
 ${\cal C}[Q(\Sigma)]\setminus\{0\}$. 
 
 Let $\{\varphi_{\lambda}|\lambda\in P^+(2\Sigma)\}$ be a $W_{\Theta}$ 
 invariant zonal spherical family associated to $\gg,\gg^{\theta}$
(See [L2], [L3], and [L4, Section 4]). Recall that 
$\varphi_{\lambda}\in {\cal C}[P(2\Sigma)]$ for each $\lambda$.   Moreover, 
elements of   
 of the ring ${\cal C}[P(\pi)]$ (which contains ${\cal C}[P(2\Sigma)]$)
 can be thought of as functions on $  U^0$ 
 where $$z^{\lambda}(\tau(\mu))=q^{(\lambda,\mu)}$$ for all 
 $\lambda\in P(\pi)$ and $\tau(\mu)\in T$.  
 
\begin{definition} ([L4, Theorem 4.1])
  Fix $B\in {\cal B}$.  The {\it radial component map} ${\cal X}$  is 
a function       ${\cal X}:  
 \check U^B\mapsto{\cal C}(Q( \Sigma)){\cal A}$ such that
 $$(\varphi_{\lambda}*{\cal 
 X}(u))(t)=z^{\lambda}({\cal 
 P}_B(u)) \varphi_{\lambda} (t)\eqno{(1.2)}$$
 for all $\lambda\in P^+(2\Sigma)$, $t\in {\cal A}$, and $u\in\check 
 U^B$.
\end{definition}
 Note that (1.2) ensures that  $\varphi_{\lambda}$ is an eigenvector for the 
 action of ${\cal X}(u)$ with eigenvalue equal to $z^{\lambda}({\cal 
 P}_B(u))$. 
  
 We have the following results from [L4] concerning the images of 
 $\check U^B$ under these two maps. (See [L4, Theorem 2.4, Theorem 2.6 and 
 Corollary 4.2].) 
 \begin{theorem} Let $B\in {\cal B}$.  The map
  ${\cal P}_B$ restricts to an algebra homomorphism from $
  \check U^B$ into ${\cal C}[{\cal A}]$ with  kernel  equal to $(B\check 
  T_{\Theta})_+\check U\cap \check U^B$.
   Moreover   there is   an algebra isomorphism from ${\cal X}(\check U^B)$ 
  onto ${\cal P}_B(\check U^B)$ given by ${\cal X}(u)\mapsto
  {\cal P}_B(u)$ for all $u\in \check U^B$.  
  \end{theorem}

 Unless  state otherwise, we assume that $B$ is an arbitrary fixed 
 subalgebra of ${\cal B}$.  At the end of the paper, we include an appendix with definitions of 
 symbols introduced later in the paper.   For more information about 
 notation and for undefined notions,   the reader is referred 
   to [Jo], [L3], and [L4].

\section{Weyl group invariance}
    
    In [L4],   ${\cal P}_B(\check U^B)$ is shown to be 
    invariant under a dotted action of the restricted Weyl group $W_{\Theta}$  whenever
    $\gg,\gg^{\theta}$ is not of type EIII, EIV, EVII, or
    EIX. The purpose of this section is to provide   criteria for the  extension of this result to the 
    remaining four cases. Recall that  the center $Z(\check U)$ of 
    $\check U$ is a 
    subalgebra of $\check U^B$ ([L3, Lemma 3.5]).   Furthermore, it is 
    shown in [L4, Lemma 5.3] that 
    ${\cal P}_B(Z(\check U))$ is $W_{\Theta}.$ invariant. Exploiting 
    the fact that ${\cal P}_B(Z(\check U))$ 
    is ``large'' inside of ${\cal P}_B(\check U^B)$, we find necessary 
    conditions in Theorem 2.5  
    which guarantee the $W_{\Theta}.$ invariance of ${\cal P}_B(\check U^B)$.
    Then, in Section 4, we show that ${\cal P}_B(\check U^B)$ satisfies 
    these conditions.   
    
   We recall some notation from [L4]. Let $\rho$ denote the half sum of 
   the positive roots in $\Delta$. The
    dotted action of $W_{\Theta}$ on ${\cal C}[{\cal A}]$ is defined 
    by $$w.q^{(\rho, 
    \lambda)}\tau{(\lambda)}=q^{(\rho,w\lambda)}\tau{(w\lambda)}$$ for all 
   $w\in W_{\Theta}$ and $\lambda\in P(\Sigma)$.     Given $\eta\in P^+(\Sigma)$, set $$\hat 
    m(2\eta)=\sum_{\gamma\in 
    W_{\Theta}\eta}q^{(\rho,2\gamma)}\tau(2\gamma).$$   Note that 
    the set $\{\hat m(2\eta)|\eta\in P^+(\Sigma)\}$ is a basis for 
    ${\cal C}[{\cal A}]^{W_{\Theta}.}$.
    
One of the important tools used to analyze the radial components in 
[L3] and [L4]  is a degree function and  related filtration (see 
[L4, (4.1), (4.2), and (4.3)]). We 
briefly recall how this degree function behaves on ${\cal C}[{\cal A}]$. 
Suppose that $\tau(\mu)\in {\cal A}$. We can write
$\mu=\sum_im_i\tilde\alpha_i$ where the $m_i$ are rational numbers.   
Set  $\deg 
\tau(\mu)=-\sum_im_i$. Let ${\rm tip}$ be the function on ${\cal 
C}[{\cal A}]$ which computes the highest degree homogeneous term. Thus
$${\rm tip}(g)=\sum_{\{\gamma|\deg\tau(\gamma)=
\deg g\}}a_{\gamma}\tau(\gamma) $$ for all Laurent polynomials 
$g=\sum_{\gamma}a_{\gamma}\tau(\gamma)$
in ${\cal C}[{\cal A}]$.
By [L4, Lemma 4.3],
$$\{{\rm tip}(X)| X\in {\cal P}_B(\check 
U^B)\}\subseteq {\rm span} \{\tau(-2\eta)|\eta\in 
P^+(\Sigma)\}.\eqno{(2.1)}$$  Setting $\deg f=0$ for all $f\in {\cal 
C}(Q(\Sigma))$ allows us to extend  the  degree function on ${\cal C}[{\cal A}]$ 
 to a function  on ${\cal C}(Q(\Sigma)){\cal A}$.
It follows 
from [L4, (4.4)] that $${\cal X}(a)={\rm tip}({\cal P}_B(a))f+ {\rm lower\ 
degree\ 
terms}\eqno{(2.2)}$$ 
for some $f\in {\cal C}(Q(\Sigma))$.

Let ${\cal C}[[Q(\Sigma)]]$ 
denote the power series ring  consisting 
of  possibly infinite sums of the form
$\sum_{\gamma\in Q^+(\Sigma)} a_{-\gamma}z^{-\gamma}$  for
$a_{\gamma}\in {\cal C}$.   
The next lemma provides another connection between the radial 
components and the Harish-Chandra projection of $B$ invariant elements.

\begin{lemma}  For all  $u\in \check U^B$, $${\cal X}(u)\in {\cal 
P}_B(u)+\sum_{\alpha_i\in\pi\setminus\pi_{\Theta}}z^{-\tilde\alpha_i}{\cal 
C}[[Q(\Sigma)]]{\cal A}.$$
\end{lemma}

\noindent
{\bf Proof}: 
Suppose $u\in \check U^B$.    As in the proof of [L3, Theorem 3.4], we can write 
$${\cal X}(u)=\sum_{i=1}^m\sum_{\gamma\leq\beta_i} a_{\gamma}z^{\gamma}$$  where $\gamma$ 
and the $\beta_i$, for $1\leq i\leq m$,  are elements of $Q(\Sigma)$
and the $a_{\gamma}$ are Laurent polynomials in ${\cal 
C}[{\cal A}]$.  Moreover, we may assume that the $\beta_i$, $1\leq 
i\leq m$, are not 
comparable via the standard partial ordering on $\gh^*$. 

    Since  $\{a_{\beta_i}|\ 1\leq i\leq m\}$ is a finite set of   Laurent polynomials in ${\cal 
C}[{\cal A}]$, it follows that there exists $\lambda\in P^+(2\Sigma)$ 
such that $z^{\lambda}(a_{\beta_i})\neq 0$ for each $i$.  By
[L2, Lemma 4.1]  the zonal spherical function 
$\varphi_{\lambda}$ is contained in $ z^{\lambda}+\sum_{\gamma<\lambda}{\cal 
C}z^{\gamma}$.   Thus $$\eqalign{\varphi_{\lambda}*{\cal X}(u)=&
\varphi_{\lambda}*\sum_{i=1}^m\sum_{\gamma\leq\beta_i} 
a_{\gamma}z^{\gamma}\cr
\in&\sum_{i=1}^m(z^{\lambda}(a_{\beta_i}))z^{\lambda+\beta_i}+\sum_{\gamma< 
\beta_i}{\cal C}z^{\gamma}.\cr}$$  The fact that  $\varphi_{\lambda}$ is an 
eigenvector for the action of ${\cal X}(u)$ forces $m=1$,
$\beta_1=0$, and $z^{\beta_1}=1$.   The lemma now follows from 
the fact  that  for each $\lambda\in P^+(2\Sigma)$,
the eigenvalue of $\varphi_{\lambda}$ is $z^{\lambda}({\cal 
P}_B(u))$ (see Definition 1.2). $\Box$
\medskip

Given $\eta\in P^+(\Sigma)$, set $${\cal 
N}_{\eta}=\sum_{\gamma'\leq \eta}\sum_{\beta\in W\gamma'}{\cal 
C}\tau(-2\beta) {\rm \ and \ }{\cal 
N}^+_{\eta}=\sum_{\gamma'< \eta}\sum_{\beta\in W\gamma'}{\cal 
C}\tau(-2\beta).\eqno{(2.3)}$$
Consider an  element  $a\in \check U^B$ such that 
$\tip(a)=\tau(-2\eta)$ for some $\eta\in P^+(\Sigma)$.     By 
[L4, Theorem 4.1], ${\cal X}(a)$ is $W_{\Theta}$ invariant.   Hence
Lemma 2.1 implies  that  
 ${\cal 
P}_B(a)\in {\cal N}_{\eta}$.

Let $w'_0$ denote the 
   longest element of the Weyl group $W_{\Theta}$.  To make the 
   notation of the rest of this section easier to read, we assume 
   that $w'_0=-1$.   It should be noted that this is true for the four 
   types of irreducible symmetric pairs which are the primary interest 
   of this section.  Moreover, the results are easily extended to the 
   case when $w_0'$ is not equal to $-1$. 
   
Suppose that $\mu\in P^+(\pi)$.   Recall that there exists a unique   central 
   element $z_{2\mu}$ in $\tau(2\mu)+(\adr U_+)\tau(2\mu)$ (see [L4, 
   Section 5]).  By [L4, Lemma 5.1] and subsequent discussion, we 
 have $${\cal P}_B(z_{2\mu})\in\hat m(2\tilde\mu)+{\cal 
 N}_{\tilde\mu}\eqno{(2.4)} $$ up to a nonzero scalar. Note 
 that ${\rm tip}(z_{2\mu})={\rm tip\ }\hat m(2\tilde\mu)= 
q^{(\rho,-2\tilde\mu)}\tau(-2\tilde\mu)$. Moreover, by [L4, Lemma 
 5.3],  ${\cal P}_B(z_{2\mu})$ is invariant under the dotted action 
 of $W_{\Theta}$.

\begin{lemma}  If $a\in \check U^B$ such that  $\tip({\cal 
P}_B(a))=q^{(\rho,-2\eta)}\tau(-2\eta)$ then
$${\cal P}_B(a)=\hat{m}(2\eta))+{\cal 
N}^+_{\eta}.$$  Moreover, if $\tip({\cal 
P}_B(a))=q^{(\rho,-2\eta)}\tau(-2\eta)$ where $\eta$ is a minuscule or pseudominuscule
restricted fundamental weight, then ${\cal P}_B(a)\in {\cal C}[{\cal 
A}]^{W_{\Theta}.}$
    \end{lemma}
    
 \noindent
 {\bf Proof:} 
Assume that $u\in \check U^B$ so that $\tip(u)=q^{(\rho,-2\eta)}\tau(-2\eta)$. 
By the discussion preceding the lemma, ${\cal P}_B(u)\in{\cal N}_{\eta}$.
Thus, there exist scalars $a_{\beta}$ such that  $${\cal P}_B(u)\in\sum_{\beta\in 
W_{\Theta}\eta}a_{\beta}\tau(-2\beta)+{\cal N}^+_{\eta}.\eqno{(2.5)}$$   The fact 
that $\tip(u)=q^{(\rho,-2\eta)}\tau(-2\eta)$ ensures that $a_{\eta}$ is 
nonzero. Assume that  $\sum_{\beta\in 
W_{\Theta}\eta}a_{\beta}\tau(-2\beta)$ is not $W_{\Theta}.$ invariant.

By [L4, Theorem 3.7],  there 
exists $\gamma$ and $\lambda$ in $P^+(\pi)$ such that 
$\tilde\lambda=\tilde\gamma+\eta$.   Since ${\cal P}_B$ is an 
algebra map (Theorem 1.3), we have ${\cal P}_B(uz_{2\gamma})={\cal 
P}_B(u){\cal P}_B(z_{2\gamma})$.  Hence, by (2.4) and (2.5) we can write
$${\cal P}_B(uz_{2\gamma})=\sum_{\beta\in 
W_{\Theta}\tilde\lambda}a_{\beta}'\tau(-2\beta)+{\cal 
N}_{\tilde\lambda}^+.\eqno{(2.6)}$$  Furthermore, our assumptions on the 
$a_{\beta}$  ensure  that 
  $\sum_{\beta\in 
W_{\Theta}\tilde\lambda}a_{\beta}'\tau(-2\beta)$ is   not 
$W_{\Theta}.$ invariant. 

Set $X= {\cal 
P}_B(uz_{2\gamma})-a'_{\tilde\lambda}q^{(\rho,2\tilde\lambda)}{\cal P}_B(z_{2\lambda})$.
By (2.4) and (2.6), we have $$X\in 
\sum_{\beta\in 
W_{\Theta}\tilde\lambda, \beta\neq \tilde\lambda}b_{\beta}\tau(-2\beta)+{\cal 
N}^+_{\tilde\lambda}\eqno{(2.7)}$$ 
for some set of scalars $b_{\beta}$ which are not all equal to zero.
By [L4, Lemma 4.3], $\tip (X)$ is a linear 
combination of elements in the set $\{\tau(-2\gamma)|\gamma\in 
P^+(\Sigma)\}$.   Since $\lambda\in P^+(\pi)$, it follows from [L4, 
Lemma 3.3 ] that  $\tilde\lambda\in P^+(\Sigma)$.  Hence if   $\beta\in 
W_{\Theta}\tilde\lambda\setminus\{\tilde\lambda\}$, then $\beta$  is not an element of 
$P^+(\Sigma)$. Thus $\tip (X)$ is contained in
${\cal N}_{\tilde\lambda}^+$.  
It now follows from  the  $W_{\Theta}$ invariance  of ${\cal X}(\check 
U^B)$ ([L4, Theorem 4.1]), (2.2), Lemma 2.1, and (2.7)   that 
${\cal 
P}_B(X)\in 
{\cal N}^+_{\tilde\lambda}.$  However, $\tau(-2\beta)\notin{\cal 
N}_{\tilde\lambda}^+$ 
for any $\beta$ in $W_{\Theta}\tilde\lambda$.
This  contradiction proves the first assertion of the lemma.

Assume now that  $\eta$ is either a minuscule or pseudominuscule 
restricted fundamental weight
(see [M] or [L3, Section 7]).
It follows that   the set $\{\gamma\in 
P^+(\Sigma)|\gamma<\eta\}$ is a subset of $\{0\}$. Thus
${\cal N}^+_{\eta}$ is a subset of ${\cal C}$.  Therefore ${\cal P}_B(u)$ 
equals 
  $\hat m(2\eta)+c$ for some scalar $c$.
$\Box$

\medskip
Now suppose that $a\in \check U^B$ and $${\rm tip}\ ({\cal 
P}_B(a))=\sum_ja_jq^{(\rho, -2\gamma_j)}\tau( -2\gamma_j)$$ with 
$\gamma_j\in P^+(\Sigma)$.   It follows 
that each $\tau(\gamma_j)$ must have the same degree as ${\cal 
P}_B(a)$.  This implies 
that the $\gamma_j$ are pairwise incomparable with respect to the 
standard partial order on $P(\Sigma)$.   The proof of Lemma 2.2 can be 
easily generalized to show that  
$${\cal P}_B(a)=\sum_ja_j\hat  m(2\gamma_j))+\sum_j{\cal 
N}_{\gamma_{j}}^+.$$

It is sometimes useful to consider the ``opposite'' degree function 
on ${\cal C}[{\cal A}]$.   In particular, set ${\rm odeg} 
\tau(\eta)=m$ where $\eta=\sum_im_i\tilde\alpha_i$ and $m=\sum_im_i$.
We define a function ``top'' on ${\cal C}[{\cal A}]$ in a manner 
analogous to the tip function.   
Given $g=\sum_{\gamma}a_{\gamma}\tau(\gamma)$ in ${\cal C}[{\cal A}]$, 
set $${\rm top}(g)=\sum_{\{\gamma|{\rm odeg\ }\tau(\gamma)={\rm odeg\ } 
g\}}a_{\gamma}\tau(\gamma).$$  Recall that ${\rm tip\ }\hat m(2\eta)= 
q^{(\rho,-2\eta)}\tau(-2\eta)$. On the other hand, ${\rm top\ }\hat 
m(2\eta)=q^{(\rho,2\eta)}\tau(2\eta)$.   
The following description of 
the set of ``tops'' of the images of elements in $\check U^B$ 
under ${\cal P}_B$ is analogous to the description of the set of 
``tips'' given in [L4, Lemma 4.3].

\begin{lemma}
    The set $\{{\rm top}({\cal P}_B(a)|a\in \check U^B\}$ is a subset 
    of ${\rm span}\{\tau(2\gamma)|\gamma\in P^+(\Sigma)\}.$
    Moreover, ${\rm top}({\cal P}_B(a)) =w'_0.({\rm tip}({\cal 
    P}_B(a))$ for all $a\in \check U^B$.
    
   \end{lemma}

\noindent
{\bf Proof:} Suppose that $a\in \check U^B$.  There exist elements 
$\gamma_1,\dots, \gamma_s$ in $P^+(\Sigma)$ and scalars $a_1,\dots, 
a_j$ such that  $${\rm tip}\ ({\cal 
P}_B(a))=\sum_{j=1}^sa_jq^{(\rho,-2\gamma_j)}\tau(-2\gamma_j).$$    
By  the discussion preceding the lemma,
$${\cal P}_B(a)=\sum_{j=1}^sa_j\hat  m(2\gamma_j))+\sum_{j=1}^s{\cal 
N}_{\gamma_{j}}^+.$$
Suppose that $\deg {\cal P}_B(a)=m$.  It follows that $\deg\tau(-2\gamma_j)=m$ 
for each $j$.   Hence ${\rm odeg}\ \tau(2\gamma_j)=m$ for 
all $1\leq j\leq s$.   Note that  $\gamma'<\gamma$ implies that ${\rm 
odeg}\tau(\gamma')<{\rm odeg}\tau(\gamma)$ for all $\gamma$ and 
$\gamma'$ in $P(\Sigma)$.  Thus
${\rm odeg\ } m(2\gamma_j)>{\rm odeg\ } c$ for all $c\in  {\cal 
N}_{\gamma_{j}}^+$ and all $j$.  Hence
$m={\rm odeg}({\cal P}_B(a))=\deg ({\cal P}_B(a))$. Since ${\rm 
top}(\hat m(2\eta))=q^{(\rho,2\eta)}\tau(2\eta)$ for each $\eta\in 
P^+(\Sigma)$, it follows that 
$$\eqalign{{\rm top}\ ({\cal  
P}_B(a))=&\sum_ja_jq^{(\rho, 2\gamma_j)}\tau( 2 \gamma_j)\cr=&
w_0'.(\sum_ja_j(q^{(\rho, -2\gamma_j)}\tau( -2 \gamma_j)).\Box\cr}$$

 \medskip

By Lemma 2.3 and (2.1) we have that   $$\{{\rm top}(X)| X\in {\cal P}_B(\check 
U^B)\}\subseteq {\rm span}\{\tau(2\eta)|\eta\in 
P^+(\Sigma)\}.\eqno{(2.8)}$$ In the Section 4, we show that the above 
inclusion is actually an equality.     First, we show that   equality 
in (2.8) is a sufficient 
condition for    ${\cal P}_B(\check U^B)$ to be invariant under the dotted
$W_{\Theta}$ action. The previous lemma  is all that is necessary for 
the three types of irreducible symmetric pairs:   EIII, EIV, and 
EVII.  However, the situation is more delicate when    $\gg,\gg^{\theta}$ is of  EIX. 
The following technical result handles this one special case.

\begin{lemma}   
Assume that $\gg,\gg^{\theta}$ is of type  EIX.  
Assume further   that  $$\{{\rm top}(X)| X\in {\cal P}_B(\check 
U^B)\}={\rm span}\{\tau(2\eta)|\eta\in 
P^+(\Sigma)\}.$$   Then ${\cal P}_B(\check U^B)={\cal 
C}[{\cal A}]^{W_{\Theta}.}$.
\end{lemma}

\noindent
{\bf Proof:}
Checking the list of irreducible symmetric pairs in [A], we see that 
$\gg$ is of type E8 and  $\Sigma$ is of type F4. Note that  ${\cal N}_{\eta}=\sum_{\gamma'\leq \eta}\sum_{\beta\in W\gamma'}{\cal 
C}\tau(2\beta)$ and 
${\cal N}^+_{\eta}=\sum_{\gamma'< \eta}\sum_{\beta\in W\gamma'}{\cal 
C}\tau(2\beta)$ for all $\eta\in P^+(\Sigma)$.
 Write $\omega_1,\dots, \omega_8$ for the fundamental weights 
for the root system of $\gg$.   Write 
$\eta_1,\eta_2,\eta_3, \eta_4$ for the fundamental weights of F4.   Here we 
are assuming that both sequences are ordered using the ordering in 
[H, Chapter III]. 

  Since  ${\cal P}_B$ 
is an algebra homomorphism ([L4, Theorem 2.4]), we have 
$${\rm  top}({\cal P}_B(ab))={\rm top}({\cal 
P}_B(a)){\rm top}({\cal P}_B(b)
)$$ for all $a$ and $b$ in $\check U^B$. By the assumptions of the 
lemma, there exist  $u_1,u_2,u_3,$ and $u_4$ in $\check 
U^B$ such 
that ${\rm top}({\cal P}_B(u_i))=\tau(2\eta_i)$ for $1\leq i\leq 4$.  It follows that ${\cal 
P}_B(u_1), {\cal P}_B(u_2), {\cal P}_B(u_3), $ and ${\cal P}_B(u_4)$ 
generates ${\cal P}_B(\check U^B)$. Thus it is sufficient to show that   ${\cal 
P}_B(u_i)\in {\cal 
C}[P(\Sigma)]^{W_{\Theta}.}$ for   $1\leq i\leq 4$.

  It is straightforward to check  that $\tilde\omega_1=2\eta_4$, 
$\tilde\omega_2=\eta_4+\eta_3$,
$\tilde\omega_3=\eta_4+2\eta_3,\tilde\omega_4=2\eta_4+2\eta_3,$
$\tilde\omega_5=2\eta_4+\eta_3$,
$\tilde\omega_6=2\eta_3$, $\tilde\omega_7=\eta_2$, and 
$\tilde\omega_8=\eta_1$.  Set $c_1=z_{2\omega_8}$ and 
$c_2=z_{2\omega_7}$.  By (2.4),
 ${\rm top}({\cal P}_B(c_i))=q^{(\rho,2\eta_i)}\tau(2\eta_i)$  for 
 $i=1,2$ and by [L4, Lemma 5.3], both ${\cal P}_B(c_1)$ and ${\cal 
 P}_B(c_2)$ are $W_{\Theta}.$ invariant.
Note that $\eta_4$ is pseudominuscule. Hence, Lemma 2.2 ensures that  there exists $c_4$ 
in $\check U^B$ such that ${\cal P}_B(c_4)=\hat 
m(2\eta_4)$.  In particular, ${\cal P}_B(c_4)$ is $W_{\Theta}.$ 
invariant and ${\rm top}({\cal P}_B(c_4))=q^{(\rho,2\eta_4)}\tau(2\eta_4) $.

Now choose $c_3\in \check U^B$ such that ${\rm top}({\cal 
P}_B(c_3))=q^{(\rho,2\eta_3)}\tau(2\eta_3))$.    Note that the only nonzero dominant integral 
weights less than $\eta_3$ are   $\eta_1$  and $\eta_4$.  Moreover, 
$\eta_4<\eta_1$. It follows that ${\cal N}_{\eta_3}^+={\cal 
N}_{\eta_1}$.  Thus by Lemma 2.2, ${\cal P}_B(c_3)\in\hat 
m(2\eta_3)+{\cal N}_{\eta_1}$. Furthermore, there exist scalars 
$a_{\beta},$ for $ \beta\in W_{\Theta}\eta_1$, such that 
$${\cal P}_B(c_3)\in \hat m(2\eta_3)+\sum_{\beta\in 
W_{\Theta}\eta_1}a_{\beta}\tau(2\beta)+{\cal N}_{\eta_4}.\eqno{(2.9)}$$

 Assume  
 that $\sum_{\beta\in 
W_{\Theta}\eta_1}a_{\beta}\tau( 2\beta)$ is not $W_{\Theta}.$ invariant.
   Consider the element $c_3c_4$ of $\check U^B$.   Note  that $\eta_2$ is the only    dominant integral 
weight less than 
$\eta_3+\eta_4$ and greater than $\eta_1+\eta_4$.  Since ${\cal P}_B$ 
is an algebra homomorphism, we have that ${\cal P}_B(c_3c_4)={\cal P}_B(c_3){\cal 
P}_B(c_4)$. It follows that 
$$\eqalign{ {\cal 
P}_B(c_3c_4)\in&\ {\rm span}\{\hat 
m( 2 \eta_3+2\eta_4),\hat m( 2\eta_2)\}\cr+& \sum_{\beta\in W_{\Theta} 
(\eta_4+\eta_1)}b_{\beta}\tau(2\beta)+ {\cal N}^+_{\eta_4+\eta_1}.\cr}$$   
Furthermore, the assumption on the $a_{\beta}$ ensures that $\sum_{\beta\in W_{\Theta} 
\eta_1+\eta_4\}}b_{\beta}\tau(2\beta)$ is not $W_{\Theta}.$ invariant. 
Now $\tilde\omega_2=\eta_3+\eta_4$ and $\tilde\omega_7=\eta_2$. By 
(2.4)  there exists 
$X\in Z(\check U)$   such that 
$${\cal 
P}_B(c_3c_4+X )\in  \sum_{\beta\in W_{\Theta} 
(\eta_1+\eta_4)}b'_{\beta}\tau(2\beta)+{\cal N}^+_{\eta_4+\eta_1} $$
and $\sum_{\beta\in W_{\Theta} 
(\eta_1+\eta_4)}b_{\beta}'\tau(2\beta)$ is not $W_{\Theta}.$ invariant. This 
contradicts   Lemma 2.2.  Hence 
$${\cal P}_B(c_3)\in \hat 
m(2\eta_3)+a \hat m(2\eta_1)+{\cal N}_{\eta_4}$$ for some 
scalar $a$.

Adding a constant to $c_3$ if necessary, we may assume that $${\cal P}_B(c_3)= \hat 
m(2\eta_3)+a \hat m(2\eta_1)+ \sum_{\beta\in 
W_{\Theta}{\eta_4}}e_{\beta}\tau(2\beta)$$ for some scalars $e_{\beta}$.
  Assume further    that $\sum_{\beta\in 
W_{\Theta}\eta_1}e_{\beta}\tau(2\beta)$ is not $W_{\Theta}.$ 
invariant.  
Examining the dominant integral weights of F4 less than or equal to 
$\eta_3+\eta_4$ yields $$\eqalign{{\cal 
P}_B(c_3c_4)&\in {\cal C}\hat 
m(2\eta_3+2\eta_4)+{\cal C}\hat m(2\eta_1+2\eta_4)+{\cal C}\hat 
m(2\eta_2)\cr&+ \sum_{\beta\in W_{\Theta} 
(2\eta_4)}e'_{\beta}\tau(2\beta)+{\cal N}_{\eta_3}.\cr}$$  Moreover, our assumptions on $c_4$ and $c_3$ 
ensure that $\sum_{\beta\in W_{\Theta} 
(2\eta_1)}e'_{\beta}\tau(2\beta)$ is not $W_{\Theta}.$ invariant.   
 Note that ${\rm top}({\cal P}_B(c_1c_4))$ is  a 
 scalar multiple of $\tau(2\eta_1+2\eta_4)$
 and 
  ${\cal P}_B(c_1c_4)$ is $W_{\Theta}.$ invariant. On the 
other hand, both $\eta_3+\eta_4$ and $\eta_2$ are in 
$\widetilde{P^+(\pi)}$. Hence using (2.4) we can find
   $X\in \check U^B$ and scalars $e''_{\beta}$ such that 
$${\cal 
P}_B(c_3c_4+X )\in  \sum_{\beta\in W_{\Theta} 
2\eta_4}e''_{\beta}\tau(2\beta)+{\cal N}_{\eta_3}$$
and $\sum_{\beta\in W_{\Theta} 
\eta_1\}}e''_{\beta}\tau(2\beta)$ is not $W_{\Theta}.$ invariant  Once again 
this contradicts the previous lemma.  It follows that ${\cal 
P}_B(c_3)$ is $W_{\Theta}.$ invariant. 
$\Box$

\medskip
The next result provides the essential criterion which will be used to 
show that ${\cal P}_B(\check U^B)$ is $W_{\Theta}.$ invariant.

\begin{theorem}Assume that $\gg,\gg^{\theta}$ is of type EIII, EIV, 
EVII, and   EIX.  
Assume further   that  $$\{{\rm top}(X)| X\in {\cal P}_B(\check 
U^B)\}={\rm span}\{\tau(2\eta)|\eta\in 
P^+(\Sigma)\}.$$   Then ${\cal P}_B(\check U^B)={\cal 
C}[{\cal A}]^{W_{\Theta}.}$.
\end{theorem}

\noindent
{\bf Proof:}   Write 
$\eta_1,\dots, \eta_t$ for the fundamental weights of the restricted 
root system $\Sigma$.   
Since ${\cal P}_B$ is an algebra homomorphism ([L4, Theorem 2.4]), 
it is sufficient to find $u_1,
\dots,  u_t$ in $\check U^B$ such that ${\rm top}({\cal 
P}_B(u_i)=\tau(2\eta_i)$ and   ${\cal 
P}_B(u_i)\in {\cal 
C}[{\cal A}]^{W_{\Theta}.}$ for $1\leq i\leq t$.  This is exactly 
what is done for EIX in the previous lemma.

If $\eta_i\in \widetilde{P^+(\pi)}$, then by (2.3) there exists $u_i\in Z(\check 
U)$ such that ${\rm top}({\cal 
P}_B(u_i))=\tau(2\eta_i)$.  Moreover, ${\cal 
P}_B(u_i)\in {\cal C}[{\cal A}]^{W_{\Theta}.}$.
Now suppose that $\gg,\gg^{\theta}$ is of type EIV or EVII. A straightforward calculation shows that the only 
fundamental restricted weights not contained in $  \widetilde{P^+(\pi)}$ are 
$\eta_1$ and $\eta_2$.  Furthermore, both $\eta_1$ and $\eta_2$ 
are either  minuscule or pseudominuscule. 
  On the other hand suppose that $\gg,\gg^{\theta}$ 
is of type EIII.  Then $\Sigma$ is of type B$_2$
and   the only fundamental restricted weight not contained in    
$  \widetilde{P^+(\pi)}$ is $\eta_1$.  It is straightforward to check 
that $\eta_1$  is a pseudominuscule 
weight. The theorem now follows from Lemma 2.2 and Lemma 2.3. $\Box$

\section{An inverse to  the Harish-Chandra map}

Recall that the  ordinary Harish-Chandra map is a projection of 
$\check U$ onto $\check U^0$.   Furthermore, this 
Harish-Chandra map induces an isomorphism between 
   $Z(\check 
U)$ and a particular   subring of $\check U^0$.
There are  a number of different proofs of this result in the 
literature (see [B]).   The approach taken by [JL] (see also [Jo, 
Section 7]) establishes the  surjectivity of this isomorphism by 
using a map which lifts certain elements of $T$  to 
$Z(\check U)$. In this section,
we study an analog ${\cal L}$ of this map   associated to quantum symmetric pairs. 
Some of the material presented here is a generalization of parts of 
[L3, Section 7].

 Let $\phi$ be the Hopf algebra automorphism of $U$ which fixes 
elements in $T$ such that 
$$\phi(x_i)=q^{(-2\rho,\tilde\alpha_i)}x_i\qquad
\phi(y_i)=q^{(2\rho,\tilde\alpha_i)}y_i\eqno{(3.1)}$$ 
for all $1\leq i\leq n$.     Note that $\tilde\alpha_i=0$ whenever 
$\Theta(\alpha_i)=\alpha_i$.   Thus $\phi$ acts as the identity on 
the subalgebra ${\cal M}$ of $U$.

Set $T_{\Theta}=\{\tau(\beta)|\Theta(\beta)=\beta$ 
and $\beta\in Q(\pi)\}$.  Recall that $B$ is generated by ${\cal M}$, $T_{\Theta}$ 
and elements of the form $y_it_i+d_i\tilde\theta(y_i)t_i+ s_i  t_i$ for 
$\alpha_i\in \pi\setminus \pi_{\Theta}$ and suitably chosen scalars 
$d_i$ and $s_i$ in ${\cal R}$ ([L4, Section 1]).  Here $\tilde\theta$ is a lift of the involution 
$\theta$ to a ${\bf C}$ algebra automorphism of $U$.  Moreover, ${\cal 
M}$ and $T_{\Theta}$ are subsets of each coideal subalgebra in ${\cal B}$.
(See [L2, Section 7] and [L4, 
Section 1] for more details.)  

  Let $\sigma$ denote the antipode and ${\it \Delta}$ the 
comultiplication of $U$ (see [Jo, 3.2.9] or [L1, Section 7]).  The following lemma generalizes 
[L3, Lemma 7.3].

\begin{lemma} Suppose that $a\in \check U$ and $b\in B_+$.
    Then $(\adr b)a$ is contained in 
$\phi(B_+) \check U +\check UB_+$.
\end{lemma}

\noindent
{\bf Proof:}  
Note that  ${\cal M}T_{\Theta}$ is a Hopf subalgebra of 
$U$ and  a subalgebra of both $\phi(B)\check T_{\Theta}$ and 
$B\check T_{\Theta}$.  Hence $$(\adr ({\cal M}T_{\Theta})_+)a\subset 
 {\cal M}T_{\Theta}a {\cal M}T_{\Theta}\subset 
 (\phi(B))_+a+\check UB_+.$$ 

 Set $C_i=y_it_i+d_i\tilde\theta(y_i)t_i+ s_i(t_i-1)$  for 
 $\alpha_i\in \pi\setminus \pi_{\Theta}$ where the $s_i$ and $d_i$ are 
 chosen so that $C_i\in B$.  
 Note that  $C_i$ is in $B_+$ for each $i$.  Now
$\sigma(s_i(t_i-1))=s_i(t_i^{-1}-1)$. Hence the proof of  [L3, Lemma 7.3] 
shows that  $$\eqalign{\sigma(C_i)&\in 
-q^{(2\rho,\alpha_i)}y_i-d_iq^{(2\rho,\Theta(\alpha_i))}\tilde\theta(y_i)+s_i(t_i^{-1}-1)+
({\cal M}T_{\Theta})_+ U\cr
&= q^{(\rho, \Theta(\alpha_i)+\alpha_i)}\phi(C_i)t_i^{-1}+({\cal 
M}T_{\Theta})_+U.\cr}$$ Furthermore, by the arguments in [L3, Lemma 7.3] 
we have that $$\eqalign{{\it\Delta}(y_it_i+d_i\tilde\theta(y_i)t_i)&\in 
(y_it_i+d_i\tilde\theta(y_i)t_i)\otimes 1+t_i\otimes 
(y_it_i+d_i\tilde\theta(y_i)t_i)\cr&+U\otimes ({\cal 
M}T_{\Theta})_+.\cr}$$  
Now ${\it\Delta}(t_i-1)=t_i\otimes 
t_i-1\otimes 1=(t_i-1)\otimes 1+t_i\otimes (t_i-1).$
Hence $${\it\Delta}(C_i)\in C_i\otimes 1+t_i\otimes C_i+U\otimes 
({\cal M}T_{\Theta})_+.$$ It follows that $$\eqalign{(\adr C_i)a&\in 
-\sigma(C_i)a+t_i^{-1}aC_i+
Ua({\cal M}T_{\Theta})_+\cr&\subset -q^{(\rho, \Theta(\alpha_i)+\alpha_i)}\phi(C_i)t_i^{-1}a+({\cal 
M}T_{\Theta})_+ Ua+t_i^{-1}aC_i+\check U({\cal M}T_{\Theta})_+\cr
&\subset \phi(B_+)\check U+\check UB_+\cr}$$ for all $a\in \check U$.  The lemma now follows 
from the fact that $B$ is generated by the $C_i, \alpha_i\in 
\pi\setminus \pi_{\Theta}$ and ${\cal M}T_{\Theta}$.
$\Box$

\medskip

Let $\chi$ denote the Hopf algebra automorphism of $U$ defined by 
$\chi(x_i)=q^{(\rho,\tilde\alpha_i)}x_i$ and 
$\chi(y_i)=q^{-(\rho,\tilde\alpha_i)}y_i$ for all $1\leq i\leq n$.  
Note  that $\phi=\chi^{-2}$. 
Note further that 
$(\rho,\tilde\alpha_i)=(\tilde\rho,\tilde\alpha_i)=(\tilde\rho,\alpha_i)$.
  Hence
$\tau(\tilde\rho)x_i\tau(-\tilde\rho)=\chi(x_i)$
and $\tau(\tilde\rho)y_i\tau(-\tilde\rho)=\chi(y_i)$ for all 
$1\leq i\leq n$.  It follows that 
$\tau(\tilde\rho)u\tau(-\tilde\rho)=\chi(u)$ for all $u\in U$.

\begin{lemma} We have the following equality of sets:
    $$\tau(\tilde\rho)(\phi(B_+)\check U+\check 
UB_+)= \chi^{-1}(B_+)\tau(\tilde\rho)\check U+\tau(\tilde\rho)\check 
UB_+.$$
    \end{lemma}
    
    \noindent 
     {\bf Proof:} By (3.1), we have that
    $\phi(u)=\tau(-2\tilde\rho)u\tau(2\tilde\rho)$ for all $u\in U$.   Hence
     $\tau(\tilde\rho)\phi(B)\tau(-\tilde\rho)=
     \tau(-\tilde\rho)B\tau(\tilde\rho)=\chi^{-1}(B)$.
     It follows that $\tau(\tilde\rho)(\phi(B))\check 
     U=\chi^{-1}(B)\tau(\tilde\rho)\check U$.  
   $\Box$
   
   \medskip

It is sometimes helpful to use
 a slightly different form of the 
  Harish-Chandra map associated to $B$. Let $N^-$ be the subalgebra 
  of $G^-$ generated by the set $(\ad {\cal M}\cap G^-){\cal 
  C}[y_it_i|\ \alpha_i\notin\pi_{\Theta}]$.
In particular, by [L4, Theorem 2.2 and (2.3)],  we have an inclusion
$$\check U\subset (\check U(B\check T_{\Theta})_++N^-_+\check {\cal 
A})\oplus {\cal C}[\check {\cal A}].\eqno{(3.2)} $$ 
Let 
${\cal P}'_{B }$ be the projection of $\check U$ onto 
${\cal C}[\check {\cal A}]$ using the direct sum decomposition of (3.2).

Let $L(\lambda)$ denote the finite-dimensional simple $U$ module of 
 highest weight $\lambda$ where $\lambda\in P^+(\pi)$.

\begin{lemma}
    For all $c\in \check U^B$, we have ${\cal P}_B(c)={\cal P}'_B(c)$. 
    Moreover, if $\lambda\in P^+(2\Sigma)$ and $\xi_{\lambda}$ is a nonzero $B$ invariant element 
    of $L(\lambda)$, then $$u\xi_{\lambda}=z^{\lambda}({\cal 
	 P}_B(u))\xi_{\lambda}$$ for all $u\in \check U^B$.
    \end{lemma}

    \noindent
    {\bf Proof:}  Let $\kappa$ be the 
  antiautomorphism of $U$ which      restricts to an algebra antiautomorphism of $B$
 as described in ([L2, Theorem 3.1]). Note that $\kappa(t)=t$ for 
 all $t\in \check T$. Moreover, $\kappa$ can be extended to $\check 
 U\check {\cal A}\check T_{\Theta}$ so that $\kappa(t)=t$ for all $t$ 
 in
 $  \check {\cal A}\check T_{\Theta}$.    As explained in the proof 
 of [L3, Theorem 2.2], we have $\kappa(N^+)=N^-$.  Hence
 applying $\kappa$ to (1.1) yields (3.2).   The first assertion of 
 the lemma now follows from 
 the fact that   $\kappa$ restricts to  the identity 
 on ${\cal C}[\check{\cal A}]$. 
 
Let $\xi_{\lambda}^*$ denote the $B$ invariant vector of 
$L(\lambda)^*$.   In [L3, preceding Theorem 3.6], it is shown that 
$\xi_{\lambda}^*c=z^{\lambda}({\cal P}_B(c))\xi^*_{\lambda}$.   
Applying $\kappa$  and switching the roles of $\xi_{\lambda}^*$ and 
$\xi_{\lambda}$  in the argument in [L3] yields $u\xi_{\lambda}=z^{\lambda}({\cal 
	 P}'_B(u))\xi_{\lambda}$.  The second assertion of the lemma now follows from the 
	 first.   $\Box$

\medskip

Let $F_r(\check U)$ denote the locally finite part of $\check U$ with 
respect to the right adjoint action.  Recall that the action of $(\adr B)$ on $F_r(\check U)$ is 
semisimple ([L4, Theorem 1.1]).   As explained in [L3, Section 7, 
before Lemma 7.5]), we have that  $$F_r(\check U)=\check U^B\oplus (\adr 
B_+)F_r(\check U).\eqno{(3.3)}$$ 

\begin{definition} The map ${\cal L}$ is  the projection map from 
$F_r(\check U)$ onto $\check U^B$ using the direct sum 
decomposition (3.3).  Moreover, given
$a\in F_r(\check U)$ such that  $a\notin (\adr B_+)F_r(\check U)$, 
we have that ${\cal L}(a)$ is the unique element of $\check U^B$ which is 
contained in  
$ a+(\adr B_+)a $ (see [L3, Section 7, before Lemma 7.5]).
\end{definition}

  Given $\lambda\in P^+(2\Sigma)$, let  
  $g_{\lambda}$ denote the zonal spherical function 
  $g^{\lambda}_{\chi^{-1}(B),B}$ (see [L2, Section 4] or [L4, Section 
  4]).   Recall that elements of $R_q[G]$ can be viewed as 
functions on $\check U$.  Write $\varphi_{\lambda}$ for the 
  image of $g_{\lambda}$ in ${\cal C}[P(\pi)]$ 
  obtained by restricting   to $U^0$. By [L2, 
 Corollary 5.4 and Theorem 6.3], each $\varphi_{\lambda}$ is in ${\cal 
  C}[P(2\Sigma)]$ and is invariant under the action of $W_{\Theta}$.  
  Moreover, the zonal spherical family $\{\varphi_{\lambda}|\lambda\in P^+(2\Sigma)\}$ 
  is a basis for ${\cal C}[P(2\Sigma)]^{W_{\Theta}}$.
  
Let $\xi_{\lambda}$ be a nonzero $B$ 
 invariant vector of $L(\lambda)$.
 Let $\zeta^*_{\lambda}$ be a nonzero $\chi^{-1}(B)$ invariant vector 
 of $L(\lambda)^*$.  By [L2, Section 4], rescaling 
 if necessary,
 $g_{\lambda}$ can be identified with the function on $\check U$ which sends 
 $u$ to $\zeta_{\lambda}^*(u\xi_{\lambda})$ for all $u\in \check U$.

\begin{lemma}   Suppose there exists $b\in \check U(B\check T_{\Theta})_+$
    and   $a\in  \check U^o$ such that $a+b\in  F_r(\check 
U)$.
Then 
$$z^{\lambda}({\cal P}_{B}({\cal L}(a+b))\varphi_{\lambda}(\tau(\tilde\rho))=\varphi_{\lambda}(a\tau(\tilde\rho))$$
for all $\lambda\in P^+(2\Sigma)$.
\end{lemma}

\noindent
{\bf Proof:} By Lemma 3.3 and the definition of $\varphi_{\lambda}$ we 
have that 
$$\eqalign{g_{\lambda} (\tau(\tilde\rho) {\cal L}(a+b))=&z^{\lambda}({\cal 
P}_{B}({\cal L}(a+b))g_{\lambda} (\tau(\tilde\rho))\cr=&z^{\lambda}({\cal 
P}_{B}({\cal L}(a+b))\varphi_{\lambda}(\tau(\tilde\rho)).\cr}$$  By Lemma 3.1
and Lemma 3.2,
$g_{\lambda} (\tau(\tilde\rho)c)=0$ for all $c\in 
(\adr B_+)\check U$.  Note further that $g_{\lambda}(\tau(\tilde\rho)b)=0$ since 
$\tau(\tilde\rho)b\in \check U(B\check T_{\Theta})_+$.
It follows from Definition 3.4 that 
$$g_{\lambda}(\tau(\tilde\rho){\cal L}(a+b))=
g_{\lambda}(\tau(\tilde\rho)(a+b))=
\varphi_{\lambda}(\tau(\tilde\rho) a).\Box$$

\medskip

The decomposition (3.2) is established using [L4, Lemma 2.1] which is particularly well 
 suited to computing the tip of elements. On the other hand, it is 
  easier to compute the   top  of elements constructed using 
 the right adjoint action.
 Thus it is necessary to   transform the 
 information contained in [L3, Lemma 2.1] and [L4, Lemma 2.1].
This is done in the next lemma after we introduce some notation, 
mostly from [L4].

 Recall that  $N^+$ is a subalgebra of $U^+$ and $N^-$ is a 
 subalgebra of $G^-$.  Set $\hat 
 N^+=\sum_{\gamma\in Q^+(\pi)}N^+_{\gamma}\tau(-\gamma)$
 and $\hat 
 N^-=\sum_{\gamma\in Q^+(\pi)}N^-_{-\gamma}\tau(-\gamma)$.  Let
    $\kappa$ be the antiautorphism defined in [L2, Theorem 3.1] and used in the proof of Lemma 3.3.
Recall that  $\kappa$ fixes elements of $T$ and 
     $\kappa(N^+)=N^-$.
 Hence $\kappa(\hat N^+)=\hat N^-$.

 Set $T'_{\leq}$ equal to the 
 set $\{\tau(-\gamma)|\gamma\in Q^+(\pi)$ and $\tilde\gamma\in  P(2\Sigma)\}.$
 (See the discussion preceding [L4, Lemma 2.1] for the definition of a similar 
 set.)
 Given 
$\beta\in Q^+(\pi)$ and a vector subspace $S$ of $\check U$,
we write $S_{\beta,r}$ for the sum of all 
weight spaces $S_{\beta'}$ with $\tilde\beta'=\tilde\beta$.  
Let $G^+$ denote the subalgebra of $U$ generated by $x_it_i^{-1}$, 
for $1\leq i\leq n$.  Similarly, let
  $U^-$ denote the subalgebra of $U$ generated by $y_i, 1\leq 
  i\leq n$.

\begin{lemma}  For all $\beta,\gamma\in Q^+(\pi)$ and $Y\in 
U^-_{-\beta}G^+_{\gamma}$, we have
\begin{enumerate}
    \item[(i)]$Y\in 
B\hat N^+_{\beta+\gamma,r} +\sum_{\tilde\beta'<\tilde\beta+\tilde\gamma}
  BT'_{\leq}\hat N^+_{\beta',r}.$
\item[(ii)]$Y\in \hat 
N^-_{\beta+\gamma,r}B +\sum_{\tilde\beta'<\tilde\beta+\tilde\gamma}
\hat N^-_{\beta',r}T'_{\leq}B.$
\end{enumerate}
    \end{lemma}
    
    \noindent
    {\bf Proof}  As in the proof of [L3, Lemma 2.1], we can reduce to 
    the case when $Y\in U^-_{-\beta}$.  Furthermore, we may assume 
    that $Y$ is a monomial of the form $y_{i_1}\cdots y_{i_m}$ with 
    $y_{i_1}\notin B$.  Recall that $B$ contains elements of the form 
    $B_i=y_it_i+d_i\tilde\theta(y_i)t_i+s_it_i$   
    where $s_i$ and $d_i$ are scalars for each $\alpha_i\notin 
    \pi_{\Theta}$.  Thus
    $$\eqalign{Y=&B_{i_1}t_{i_1}^{-1}y_{i_2}\cdots 
    y_{i_{m}}-d_i\tilde\theta(y_{i_2})y_{i_1}\cdots 
    y_{i_{m}}\cr-&s_it_{i_1}y_{i_2}\cdots 
    y_{i_{m}}.\cr}$$  The proof of (i)  now follows using  induction 
     along the lines of   [L3, Lemma 2.1] and [L4, Lemma 2.1].  
     Assertion (ii) follows by applying $\kappa$ to (i).
    $\Box$
    
    \medskip

    Set ${\cal A}_{\leq}=\{\tau(-\gamma)|\ \gamma\in Q^+(\Sigma)\cap 
    P(2\Sigma)\}$.  Note that $${\cal 
    P}_B(T'_{\leq})=\{\tau(-\tilde\beta)|\ \tau(\beta)\in 
    T'_{\leq}\}={\cal A}_{\leq}.\eqno{(3.4)}$$
Now consider   $\tau(\beta)\in \check T$.   Note that $\hat N^-_+$ is a subset of 
 $N^-_+U^0$.  Now suppose that $Y\in U^-G^+$.   It follows from Lemma 
 3.6 (ii), (3.4),  and [L4, (2.3)] that $$Y\tau(\beta)\in \check U((B\check 
 T_{\Theta})_++{\cal N}^+_+\check {\cal A})\oplus 
 {\cal 
 C}\tau(\tilde\beta){\cal A}_{\leq}.\eqno{(3.5)
 }$$

\begin{lemma}
   Let $\gamma\in P^+(\Sigma)$. Suppose that there exists $b\in 
   \check U(B\check 
   T_{\Theta})_+$ 
    such that $\tau(2\gamma)+b\in F_r(\check U)$.   Then
   $${\cal P}_B({\cal L}(\tau(2\gamma)+b)\in \tau(2\gamma){\cal A}_{\leq}.$$
  Moreover, there exist 
    elemets $\beta_1,\dots, \beta_m$ in $P^+(2\Sigma)$ and 
    scalars $a_{\beta_i}$ such that $\beta_i\leq \gamma$ for each 
    $1\leq i\leq m$ and   
    $${\rm top}({\cal P}_B({\cal L}(\tau(2\gamma)+b)))=
   \sum_{i=1}^ma_{\beta_i} \tau(2\beta_i). $$ 
    \end{lemma}
    
    \noindent
    {\bf Proof:} By Lemma 3.3, ${\cal 
    P}_B({\cal L}(\tau(2\gamma)+b))={\cal P}'_B({\cal L}(\tau(2\gamma)+b))$.   We 
    work with this second version of the Harish-Chandra map.  
    By Definition 3.4 there 
    exists $c\in B_+$ so that 
    ${\cal L}(\tau(2\gamma)+b)=\tau(2\gamma)+b+(\adr c)(\tau(2\gamma)+b)$.   
    By assumption, $b\in \check U(B\check T_{\Theta})_+$.
    Moreover,  Lemma 3.1 ensures that $(\adr c)b\in \check U(B\check T_{\Theta})_+$.  Hence 
   ${\cal 
    P}'_B({\cal L}(\tau(2\gamma)+b))={\cal P}'_B(\tau(2\gamma)+(\adr 
    c)\tau(2\gamma))$.
    
    Now $\tau(2\gamma)+(\adr c)\tau(2\gamma)\in (\adr U)\tau(2\gamma)$.
    Examining   the right adjoint action of the 
    generators of $\check U$ ([L2, (1.2)] and the relations of $\check U$  yield   $$(\adr 
    U)\tau(2\gamma)\in \sum_{\beta\in 
    Q^+(\pi)}  U^-G^+\tau(2\gamma-2\beta).$$ Thus, the first assertion  follows 
    from (3.5). The second assertion is now a consequence of (2.8). $\Box$

   \medskip

  In Theorem 2.5, a criterion was found concerning the set of 
  tops of elements in ${\cal P}_B(\check U^B)$ which ensures 
  surjectivity of the Harish-Chandra map.   The next result provides 
  another, more useful, criterion.  In particular, we show that the image of $\check U^B$ 
  under ${\cal P}_B$ is all of ${\cal C}[{\cal A}]^{W_{\Theta}.}$
  if $F_r(\check U)$ 
  contains a set of elements of a particular form.  In the next 
  section, we  find this desired set of elements in $F_r(\check 
  U)$.
    
  \begin{theorem}  Suppose that for each $\gamma\in P^+(\Sigma)$, 
  there exists $b_{\gamma}\in \check U(B\check T_{\Theta})_+$ such that 
  $\tau(2\gamma)+b_{\gamma}\in F_r(\check U)$.   Then 
  $$\{{\cal 
  P}_B({\cal L}(\tau(2\gamma)+b_{\gamma}))|\gamma\in P^+(\Sigma)\}{\rm 
  \ is \ a \ basis \ for\ }{\cal P}_B(\check U^B).$$  Moreover,
  ${\cal P}_B(\check U^B)={\cal C}[{\cal A}]^{W_{\Theta}.}$.
  \end{theorem}
  
  \noindent
  {\bf Proof:}  Suppose that $a$ is a linear combination 
  of elements in the set $\{{\cal L}(\tau(2\beta)+b_{\beta})|
  \beta\in P^+(\Sigma){\rm \ and \ }
   \beta\leq \gamma\}.$   Arguing as in   [L3, Lemma 
  7.5] using Lemma 3.5, we see that 
  ${\cal P}_B(a)=0$  if and only if $a=0$. 
  Thus, the set $\{{\cal 
  P}_B({\cal L}(\tau(2\gamma)+b_{\gamma}))|\gamma\in P^+(\Sigma)\}$ 
  is linearly independent over ${\cal C}$.

  By 
  the previous lemma, ${\cal P}_B({\cal L}(\tau(2\gamma)+b_{\gamma}))\in 
  \tau(2\gamma){\cal A}_{\leq}$.   
  Let 
  $$S_{\gamma}={\rm span}\{{\cal P}_B({\cal L}(\tau(2\beta)+b_{\beta}))|
  \beta\in P^+(\Sigma){\rm \ and \ }
   \beta\leq \gamma\}.$$   It 
  follows that the dimension of $S_{\gamma}$ over ${\cal C}$ is just 
  the cardinality of the set $\{\beta\in P^+(\Sigma)$ and $\beta\leq 
  \gamma\}$.  On the other hand,   Lemma 3.7  ensures that $\{{\rm top}(a)|a\in S_{\gamma}\}$ is a 
  subspace of ${\rm span}\{\tau(2\beta)|\beta\in P^+(\Sigma)$ and $\beta\leq 
  \gamma\}$.   Since the dimension of $\{{\rm top}({\cal P}_B(a))|a\in S_{\gamma}\}$
  is equal to the dimension of $S_{\gamma}$ and both are finite 
  dimensional, it follows that 
   $$\{{\rm top}(a)|a\in S_{\gamma}\}={\rm span}\{\tau(2\beta)|
   \beta\in P^+(\Sigma){\rm \ and \ }\beta\leq \gamma\}.$$  This forces  ${\rm top}({\cal 
   P}_B({\cal L}(\tau(2\beta)+b_{\beta})=\tau(2\beta)$ up to a nonzero 
   scalar for all $\beta\in P^+(\Sigma)$.  Hence (2.8) and the above 
   discussion yields $$\eqalign{\{{\rm top}({\cal P}_B(a))|a\in \check 
   U^B\}&\subseteq {\rm span}\{\tau(2\beta)|\beta\in P^+(\Sigma)\}\cr &=
   {\rm span}\{{\rm top}({\cal P}_B({\cal L}(\tau(2\beta)+b_{\beta})))|
   \beta\in P^+(\Sigma)\}\cr
  &\subseteq {\rm span}\{{\rm top}({\cal P}_B(a))|a\in \check 
   U^B\}.\cr}$$  The theorem now follows from Theorem 
   2.5.  $\Box$

\section{Surjectivity of the Harish-Chandra map}

In this section, we prove that the image of   $\check U^B$ 
under ${\cal P}_B$ is the entire
   invariant ring ${\cal C}[{\cal A}]^{W_{\Theta}.}$.  The approach 
   is to show that $F_r(\check U)$ contains a set of elements 
   $$\{\tau(2\gamma)+b_{\gamma}|\gamma\in P^+(\Sigma)\}\eqno{(4.1)}$$ where $b_{\gamma}\in 
   \check U(B\check 
   T_{\Theta})_+$ and then apply Theorem 3.8.  In particular, 
   we prove the following theorem and corollary,  the main results of 
   this paper.
   
   \begin{theorem}
       For each $\gamma\in P^+(\Sigma)$, there exists $b_{\gamma}\in 
      \check U (B\check T_{\Theta})_+$ such that
      $\tau(2\gamma)+b_{\gamma}\in F_r(\check U)$.  Moreover
      $$\{{\cal 
      P}_B({\cal L}(\tau(2\gamma)+b_{\gamma} ))|\ \gamma\in P^+(\Sigma)\}$$
      is a basis for ${\cal P}_B(\check U^B)$.
       \end{theorem}
       
     The next corollary is an immediate consequence of 
     Theorem 4.1, [L4, Corollary 4.2 and Thoerem 5.5], and Theorem 3.8.
     
     \begin{corollary}  For each irreducible symmetric pair 
     $\gg,\gg^{\theta}$ and each $B\in {\cal B}$, 
     the Harish-Chandra map ${\cal P}_B$ maps $\check U^B$ onto ${\cal 
     C}[{\cal A}]^{W_{\Theta}.}$.  The kernel of  the restriction of 
     ${\cal P}_B$   to $\check U^B$ is   $(B\check 
     T_{\Theta})_+\check U\cap \check U^B$.  Moreover, 
     ${\cal P}_B(Z(\check U))={\cal P}_B(\check U^B)$ if and only 
     if $\gg,\gg^{\theta}$ is not of type EIII, EIV, EVII, or EIX.
     \end{corollary}

    Suppose that we have found a set of elements $\{b_{\gamma}|\ 
    \gamma\in P^+(\Sigma)\}$ which satisfy the conditions of Theorem 4.1.
   Recall that the compact zonal spherical functions have been 
   identified with Macdonald polynomials for a large class of quantum 
   symmetric pairs (see for example [L3]).  It follows from Lemma 3.5 the 
   set $\{{\cal X}({\cal 
    L}(\tau(2\gamma)+b_{\gamma}))|\ \gamma\in P^+(\Sigma)\}$ is 
    precisely the family of commuting difference operators 
    associated to these Macdonald polynomials described in [K, Theorem 
    6.6].  Thus Theorem 4.1 and 
    the results of Section 3 provide a natural quantum interepretation 
    of     these difference operators.

 Recall the antiautomorphism  $\kappa$ defined in [L2, Theorem 3.1] 
 and 
used in the proof of Lemma 3.3.  
 Since $\kappa(t)=t$ for all $t\in T$, it follows that
  $\kappa((\adr 
  t)a)=(\adr t)\kappa(a)$ for all $a\in \check U$.  Moreover,  
  up to nonzero scalars, $\kappa((\adr x_i)a)=(\adr y_i)\kappa(a)$, 
  and $\kappa((\adr y_i)a)=(\adr y_i)\kappa(a)$,  for all $1\leq i\leq n$
  (see the proof of [L3, Theorem 2.2]). 
  It follows that $\kappa(F_r(\check U))=F_r(\check 
  U)$. Since $\kappa$ restricts to an antiautomorphism of $B$, we see 
 further that 
    $\kappa((B\check T_{\Theta})_+\check U)=\check 
  U(B\check T_{\Theta})_+$.
The strategy in this section is to  show that $F_r(\check U)$ contains a set of elements 
$$\{\tau(2\gamma)+b'_{\gamma}|\gamma\in P^+(\Sigma)\}\eqno{(4.2)}$$ where each 
$b'_{\gamma}\in (B\check T_{\Theta})_+\check U$.    Applying $\kappa$ 
to this set yields the desired set of the form (4.1).  The rest of 
Theorem 4.1 then follows from Theorem 3.8.

Recall that $F_r(\check 
U)\cap \check T$ is equal to $\{\tau(2\mu)|\mu\in P^+(\pi)\}$ (see [Jo, 
Section 7] and discussion preceding [L3, Lemma 7.2]). The next lemma shows that it is very easy to find $b_{\gamma}'$ when 
$\gamma\in \widetilde{P^+(\pi)}$. 

\begin{lemma}  Suppose that $\gamma\in P^+( \Sigma)$ and 
$\gamma=\tilde\beta$ for some $\beta\in P^+(\pi)$.   Then 
$$\tau(2\gamma)+b_{\gamma}'\in F_r(\check U)$$ where
$b_{\gamma}'=\tau(2\beta)-\tau(2\gamma)$.  Moreover, 
$b_{\gamma}'\in {\cal C}[\check T_{\Theta}]_+\check T$.
\end{lemma}   

\noindent
{\bf Proof:}  Since $\beta\in P^+(\pi)$, it follows that 
$\tau(2\beta)\in F_r(\check U)$.   Therefore 
$\tau(2\gamma)+b_{\gamma}'=\tau(2\beta)$ is in $F_r(\check U)$.   
Since $\gamma=\tilde\beta=(\beta-\Theta(\beta))/2$, we have
$$\eqalign{b_{\gamma}'&=\tau(2\beta)-\tau(2\gamma)=\tau(\beta+\Theta(\beta))\tau(2\gamma)-\tau(2\gamma)\cr
&=[(\tau(\beta+\Theta(\beta))-1]\tau(2\gamma).\cr}$$ The lemma follows 
from the fact that $(\tau(\beta+\Theta(\beta))-1\in {\cal C}[\check 
T_{\Theta}]_+$.$\Box$
\medskip

We break the proof of Theorem 4.1 into four cases.

\begin{enumerate}
    \item[(i)] $\gg,\gg^{\theta}$ is not of type EIII, EIV, EVII,   
    EIX, or CII(ii) and  
$\gg$ does not contain a $\theta$ invariant Lie subalgebra $\gr$ of 
rank greater than or equal to $7$ such that $\gr,\gr^{\theta}$ is of 
type AII. 
\item [(ii)] $\gg,\gg^{\theta}$  is of type CII(ii) or $\gg$ contains a $\theta$ invariant Lie subalgebra $\gr$ of 
rank greater than or equal to $7$ such that $\gr,\gr^{\theta}$ is of 
type AII.
\item [(iii)] $\gg,\gg^{\theta}$ is of type EIV, EVII, or EIX.
\item [(iv)] $\gg,\gg^{\theta}$ of of type EIII.
\end{enumerate}

  It should be noted that in the first two cases, Corollary 4.2 follows 
from [L4, Theorem 6.1].   Indeed, Cases (i) and (ii) correspond precisely to the situation 
when the image of the center $Z(\check U)$ under ${\cal P}_B$ is equal 
to the image of $\check U^B$ under ${\cal P}_B$.  On the other hand, 
Cases (iii) and 
(iv) rely on an analysis of symmetric pairs of type DI. Thus after 
completing Case (ii), we take a detour and focus on type DI pairs. 

\bigskip

\noindent
{\bf Case (i)}:    By [L4, Theorem 3.5], $\widetilde{P^+(\pi)}=
P^+(\Sigma)$ when $\gg,\gg^{\theta}$ is an irreducible symmetric pair 
satisfying the conditions of Case (i).  Thus this case follows from 
Lemma 4.3. 

\bigskip

For the remaining three cases, finding appropriate 
 elements $b_{\gamma}'$ in $(B\check T_{\Theta})_+\check U$
 so that $\tau(2\gamma)+b_{\gamma}'$ is in 
$ F_r(\check U)$ is more delicate. Set $t$ equal to the rank of  $\Sigma$.  Let $\mu_1,\dots, \mu_t$ denote the simple roots for the restricted 
 root system $\Sigma$.    Let $\eta_1,\dots, \eta_t$ denote 
the corresponding fundamental weights in $P^+(\Sigma)$. The next lemma 
reduces the work 
 to finding just a finite number of elements associated to the 
  weights $\eta_1,\dots, \eta_t$.
 
 \begin{lemma}   Suppose that for each $B\in {\cal B}$ there exists a 
 subset   $\{b_{i}^B|1\leq i\leq r\}$ of $(B\check T_{\Theta})_+\check U$  such that 
 $\{\tau(2\eta_i)+b_i^B|1\leq i\leq r\}$ is a subset of $F_r(\check 
 U)$.   
 Then for each $B\in {\cal B}$, $F_r(\check U)$ contains a subset of 
 the form (4.2).
 \end{lemma}

 \noindent
 {\bf Proof:}  Let $S$ be the subset of $\{\tau(2\gamma)|\ \gamma\in 
 P^+(\Sigma)\}$ such that $\tau(2\gamma)\in F_r(\check U)+(B\check 
 T_{\Theta})_+\check U$ for each $\tau(\gamma)\in S$ and every $B\in {\cal 
 B}$. To prove the lemma, it is sufficient to show that $S$ equals 
 $\{\tau(2\gamma)|\ \gamma\in 
 P^+(\Sigma)\}$.  By assumption, $S$ contains the subset  
 $\{\tau(2\eta_i)|\ 1\leq i\leq 
t\}$.  Hence, it is sufficient to 
   show that $S$ is multiplicatively closed.    
   
   Fix $\eta$ and $\beta$
so that $\tau(2\eta)$ and $\tau(2\beta)$ are both in $S$.   The map $u\mapsto 
 \tau(\eta)u\tau(-\eta)$ defines a Hopf algebra automorphism of 
 $\check U$, which we denote by $\psi$.   Note that elements of 
 $ U^0$ are fixed under the action of $\psi$.  It follows that $\psi$ permutes the elements 
 of ${\cal B}$.   Fix $B\in {\cal B}$ and let $B'\in {\cal B}$ be 
 chosen so that $\psi(B')=B$.    Choose $b$ in $(B\check 
 T_{\Theta})_+\check U$ and $c$ in 
 $(B'\check 
 T_{\Theta})_+\check U$ so that $\tau(2\eta)+b$ and $\tau(2\beta)+c$ 
 are elements of $F_r(\check U)$.  Note that $\psi(c)$ is in $(B\check 
 T_{\Theta})_+\check U$.
 Hence
 $$\eqalign{(\tau(2\eta)+b)(\tau(2\beta)+c)=&
 \tau(2\eta+2\beta)+b (\tau(2\beta)+c)+\tau(2\eta )c\cr
 =&\tau(2\eta+2\beta)+b (\tau(2\beta )+c)+\psi(c)
 \tau(2\eta)\cr
\in& \tau(2\eta+2\beta)+(B\check 
 T_{\Theta})_+\check U.\cr}$$
 Thus $S$ contains $\tau(2\beta)\tau(2\eta)=\tau(2\eta+2\beta)$. $\Box$
 
 \bigskip

\noindent
{\bf Case (ii):}
Assume that $\gg,\gg^{\theta}$ satisfies the conditions of Case 
(ii). A list of the possibilities for  $\gg,\gg^{\theta}$ can be found in 
 [L4, (3.11)-(3.15)].    Note that  $t\geq 4$ in all cases (and $t\geq 7$ 
 if $\gg,\gg^{\theta}$ is not of type CII(ii)).
 From the list in [L4, Section 3] and the classification in [A], we see that $\Sigma$ must be of 
  type A$_t$, B$_t$, or C$_t$.
  Note further that we may assume that 
  the restricted root system associated to the smaller symmetric pair 
  $\gr,\gr^{\theta}$ has rank $t-1$.

    The next lemma will 
  be used to  relate  
  the fundamental weights of the root system of type A$_j$ 
  generated by $\mu_1,\dots, \mu_j$ to the fundamental weights 
  associated to 
   $\Sigma$.
  
  \begin{lemma} Assume that $\phi$ is a root system of 
  type A$_m$, $C_m$ or $D_m$ with set of simple roots 
  $\{\beta_1,\dots, \beta_m\}$.   
  Let $\lambda_1,\dots, \lambda_m$ denote the corresponding fundamental weights in 
  the weight lattice of $\phi$.  Fix $j$ such that  $1\leq j<m$. Write
  $\lambda_1',\dots, 
  \lambda_j'$ for the fundamental weights in the root system of type $A_j$ 
  generated by $\beta_1,\dots, \beta_j$ (in the obvious order).  Then 
  $$\lambda'_k+
  k(j+1)^{-1}\lambda_{j+1}=\lambda_k$$ for all $1\leq k\leq j$.
      \end{lemma}
      
  \noindent
  {\bf Proof:} The assumptions on $\phi$ ensure that  
  $(\beta_i,\beta_{i+1})=-(\lambda_{i+1},\beta_{i+1})$ for all $1\leq i\leq 
 m-1$. 
   Now $(\lambda'_k,\beta_i)=\delta_{ik}(\beta_i,\beta_i)/2$ for all $1\leq i\leq 
 j$.  It follows from [H, Section 13.2, Table 1] that 
  $$\eqalign{\lambda'_{k }&=[(j-k+1)\beta_1+2(j-k+1)\beta_2
   +\dots+(k-1)(j-k+1)\beta_{j-1}\cr&+
  r(j-k+1)\beta_k+
  r(j-k)\mu_{k+1}+\dots+k\beta_j](j+1)^{-1}.\cr}\eqno{(4.3)}$$
  Hence 
  $(\lambda'_{k},\beta_{j+1})=k(j+1)^{-1}(\beta_j,\beta_{j+1})=-k(j+1)^{-1}
  (\lambda_j,\beta_{j+1})$.  Note that (4.3) also implies that  
  $(\lambda'_k,\beta_s)=0$ for all $s>j+1$.  This proves the lemma.  $\Box$ 
 \medskip 
 
 Define subsets $\pi_i$ 
  of $\pi$ for 
  $1\leq i\leq t-1$ by $\pi_i=\{\alpha_j|1\leq j\leq 2i+1\}$.  Let 
  $\gg_i$ denote the semisimple Lie subalgebra of $\gg$ generated by 
  the positive and negative root vectors associated to the simple 
  roots in $\pi_i$. Note that the root system of $\gg_i$ is of type 
  $A_{2i+1}$.    We have 
  $\gg_i\subseteq \gg_{i+1}$ for $1\leq i\leq t-2$.  Moreover, $\theta$
  restricts to an involution on $\gg_i$
  such that  the symmetric pair $\gg_i,\gg_i^{\theta}$ is of 
  type AII  and the rank of the corresponding 
  restricted root system  is $i$. Given $i$ such that $1\leq i\leq 
  t-1$, let $U_i$ denote the quantized 
  enveloping algebra of $\gg_i$ considered as a subalgebra of 
  $U_q(\gg)$.   In particular, $U_i$ is generated by $x_j, y_j, 
  t_j^{\pm 1}$ for all $j$ such that $\alpha_j\in \pi_i$.  Set 
  $U_0={\cal C}$ and $U=U_t$.  
  
 Let $\check U_i$ be the simply connected quantized enveloping 
 algebra of $\gg_i$.   Let $\check U^i$  be the subalgebra of $U$ generated by $U_i$ and 
 $\check T$ for $0\leq i\leq t$. Note that when $i=0$, this definition 
 agrees with the earlier definition of  $\check U^0$.
 Now $\check U_i\subseteq 
 \check U^i$ for all $i$ and $\check U^t=\check U$.  Write 
 $\omega_j$ for the fundamental weight corresponding to the simple 
 root $\alpha_j$ in $P^+(\pi)$.  
Note that $\tau(\omega_j)$ is in the center of $\check 
 U_i$ for each $j>2i+1$.   In particular, $\check U^i$ is isomorphic 
 to the tensor product $\check U_i\otimes_{\cal C}{\cal 
 C}[\tau(\omega_{2i+2})^{\pm 1},\dots, \tau(\omega_{n})^{\pm 1}]$ as an algebra.   
 Furthermore, $F_r(\check U^i)=F_r(\check U_i){\cal C}
 [\tau(\omega_{2i+2}^{\pm 1}
 ),\dots, \tau(\omega_{n})^{\pm 1}]$.

 Recall the map ${\cal L}$ which lifts elements of $F_r(\check U)$ to 
 elements of $\check U^B$ defined in Definition 3.4.  Given $i$ such 
 that $1\leq i<t$, we define the same type of map, denoted by 
 ${\cal L}_i$, for the quantum symmetric pair $U_i,U_i\cap B$.  In 
 particular, given $a\notin (\adr (B\cap U_i)_+)\check U_i$, we have that 
 ${\cal L}_i(a)$ is the unique element in $$a+(\adr (B\cap U_i)_+)a.$$ 
 Note that ${\cal C}
  [\tau(\omega_{2i+2}^{\pm 1}
  ),\dots, \tau(\omega_{n})^{\pm 1}]$ is contained in the center of 
  $\check U^i$.  Thus we can extend   
 ${\cal L}_i$   to $\check U_i$ by insisting that ${\cal 
 L}_i(ab)={\cal L}_i(a)b$ for all 
 $b$ in the ring $ {\cal C}
  [\tau(\omega_{2i+2}^{\pm 1}
  ),\dots, \tau(\omega_{n})^{\pm 1}]$.
  
  Let ${\cal P}_{B\cap U_j}$ denote the Harish-Chandra map on $\check U_j$ 
    with respect to the quantum symmetric pair $B\cap U_j, U_j$. 
   Note that  ${\cal P}_{B\cap U_j}(a)={\cal P}_{B}(a)$ for all 
   $a\in \check U_j$.  Moreover, ${\cal P}_B(ab)={\cal P}_{B\cap U_j}(a)b$ 
   for all $a\in U_j$ and $b$ in the center of $\check U^j$. 
  
The next lemma shows that $F_r(\check U)$ contains a set of 
 the form (4.2).   In particular, Theorem 4.1 for 
 Case (ii) symmetric pairs  is a consequence
 of the following  result. 

\begin{lemma} For each $1\leq i\leq t$, there exists $b_i\in ((B\cap U_{i-1})\check 
 T_{\Theta})_+\check U^{i-1}$   such that 
 $\{\tau(2\eta_i)+b_i|1\leq i\leq t\}$ is a subset of $F_r(\check 
 U)$.
 \end{lemma}

\noindent
{\bf Proof:}
 By  [L4, Lemma 3.5], we have that $\tilde\omega_1=\eta_1$.  
  Thus by  Lemma 4.3, there exists $b_1\in ((B\cap U_{0})\check 
 T_{\Theta})_+\check U^{0}$ such that $\tau(2\eta_1)+b_1\in 
 F_r(\check U)$.  If $\gg,\gg^{\theta}$ is not of type CII(ii), then by 
   [L4, Lemma 3.5] $\tilde\omega_n=\eta_t$.   Thus, if $\gg,\gg^{\theta}$ is not of type CII(ii), 
 then there also exists an appropriate choice for $b_t$. 
  Note that if $\gg,\gg^{\theta}$ is of type CII(ii), then 
$\Sigma$ is of type C$_t$. Set $t'=t$ if $\gg,\gg^{\theta}$ is not of type CII(ii) and set $t'=t-1$ 
otherwise.  The proof   is by induction on $t'$.

Suppose that $1\leq j<t'$. Assume that we have found 
 $b_1,\dots, b_j$ so that $b_i\in ((B\cap U_{i-1})\check 
 T_{\Theta})_+\check U^{i-1}$ for $1\leq i\leq j$ and
 $\{\tau(2\eta_i)+b_i|1\leq i\leq j\}$ is a subset of $F_r(\check 
 U)$.   It follows that $b_i\in \check U^i$ for each $i$ and 
 $\tau(2\eta_i)+b_i\in F_r(\check U^i)$ for $1\leq i\leq j$.

 Let $\eta_1',
 \dots,\eta_j'$ denote the fundamental weights in the root system of 
 type $A_j$ generated by $\mu_1,\dots, \mu_j$.
      Lemma 4.5 and the discussion preceding 
     this lemma implies that
  $${\cal L}_j(\tau(2\eta_k))={\cal L}_j(\tau(2\eta'_k))\tau(k(j+1)^{-1}2\eta_{j+1})$$
  for each $1\leq k\leq j$.
  Let 
  $W_{\Theta}'$ denote the Weyl group of the root system generated 
  by $\mu_1,\dots, \mu_j$.  Given a dominant integral weight $\beta$ in the weight 
  lattice of this root system, set $$\hat m'(2\beta)=\sum_{\gamma\in 
  W_{\Theta}'\beta}q^{(\rho,2\gamma)}\tau(2\gamma).$$  Note that 
  $\eta'_k$ is a minuscule fundamental weight in the weight lattice of 
  the root system generated by $\mu_1,\dots, \mu_j$.  
  By Lemma 2.2,  it follows that ${\cal P}_{B\cap 
  U_j}({\cal L}_j(\tau(2\eta'_k))$ is a nonzero scalar multiple of $\hat 
  m'(2\eta'_k)$.   Hence $${\cal P}_{B}({\cal L}_j(\tau(2\eta_k))=\hat 
  m'(2\eta'_k)\tau(k(j+1)^{-1}2\eta_{j+1})\eqno{(4.4)}$$ up to a nonzero scalar.  
  By [L4, Theorem 2.4] and (4.4), we have
  $${\cal P}_{B}({\cal L}_j(\tau(2\eta_k){\cal L}_j(\tau(2\eta_{j+1-k})))=\hat m'(2\eta'_k)\hat 
  m'(2\eta'_{j+1-k})\tau(2\eta_{j+1})\eqno{(4.5)}$$ up to a nonzero scalar for all $1\leq k\leq j$.
  
 Now let $\omega_{1j},\dots, \omega_{2j+1,j}$ denote the fundamental 
 weights corresponding to the simple roots in the root system 
 generated by $\pi_j$.  Since $j<t'$, we have that $\pi_j$ generates a 
 root system of type $A_{2j+1}$.   Moreover, checking the list [L4, 
 (3.11)-(3.15)], we see that $\pi$ is of type $A_n$,$C_n$, or $D_n$. 
 Thus  Lemma 4.5 implies that 
 $$\omega_{kj}+k(2j+2)^{-1}\omega_{2j+2}=\omega_j\eqno{(4.6)}$$ for 
 $1\leq k\leq 2j+1$. By [L4, Lemma 3.1 and Lemma 3.5] , we further have that 
 $\tilde\omega_{2k}=2\eta_k$ for all $1\leq k\leq t'$.  
 Recall ([Jo, Chapter 7] and [L4, Section 5]) that there exists a unique central element
 $z'_{j+1}$  of $\check U_j$ contained in $\tau(2\omega_{j+1,j})+
  (\adr  (U_j)_+)\tau(2\omega_{j+1,j})$. Set 
  $z_{j+1}=z_{j+1}'\tau(  \omega_{2j+2})$ It follows that 
  $z_{j+1}$ is in $(\adr U_j)\tau(2\omega_j)$ and is  central in 
  $\check U^j$.  Moreover, $${\cal P}_B(z_{j+1})={\cal P}_{B\cap \check 
  U_j}(z_{j+1}')\tau( \tilde\omega_{2j+2})={\cal P}_{B}(z_{j+1}')\tau( 
  2\eta_{j+1}).\eqno{(4.7)}$$
  
  It follows from [L4, Lemma 6.7] that  $1$ is in the   ${\cal C}$ 
  span of the set $$\{{\cal P}_{B\cap \check U_j}(z'_{j+1})\}\cup\{\hat m'(2\eta'_k)\hat 
  m'(2\eta'_{j+1-k}) |\ 1\leq k\leq 
  j\}.$$  In particular, there exists $c$ in $\check U_j^{(B\cap U_j)}$ such 
  that $c$ is a linear combination of elements in the set 
 $$\{z'_{j+1}\}\cup\{{\cal L}_j(2\eta'_k){\cal L}_j(2\eta'_{j+1-k})  |\ 1\leq k\leq 
  j\} $$ and ${\cal P}_{(B\cap U_j)}(c)=1$.  By [L4, Corollary 4.2], 
  we see that $$c-1\in (B\check T_{\Theta}\cap \check U_j)_+\check 
  U_i.\eqno{(4.8)}$$  
  
  Now consider the element $c\tau(2\eta_{j+1})$. It follows from the 
  choice of $c$, (4.5), and (4.7), that
  $\tau(2\eta_{j+1})$  is a linear combination of elements in the set
  $$\{{\cal P}_{B}(z_{j+1})\}\cup\{\hat m'(2\eta'_k)\hat 
  m'(2\eta'_{j+1-k})\tau(2\eta_{j+1})|\ 1\leq k\leq j\}.$$
  Moreover, (4.8) implies that 
    $$c\in \tau(2\eta_{j+1})+((B\cap  U_j)\check T_{\Theta})_+\check U.\Box$$

 \bigskip
  
  \noindent
  {\bf Symmetric pairs of type DI:}
     The next few lemmas will provide a method for finding 
     appropriate $B$ invariant elements inside of $F_r(\check U)$ when 
     $\gg,\gg^{\theta}$ is of type EIII, EIV, EVII, or EIX.   Note 
     that in each of these cases, $\gg$ contains a semisimple Lie 
     subalgebra $\gr$  such that $\gr, \gr^{\theta}$ is of 
     type $DI$, case (i). (In fact, for the latter three cases, $\gg$   contains more than one such 
     Lie subalgebra.)   Thus in this subsection, we analyze symmetric pairs of type 
     $DI$, case (i). This information is then pulled   back to the above 
     four symmetric pair types to finish the proof of Theorem 4.1.
     
     Let ${\cal P}$ denote the ordinary Harish-Chandra projection map 
     of $\check U$ onto $\check U^0$ using the direct sum 
     decomposition [L4, (5.1)]. Let $W$ denote the Weyl group of the 
     root system $\Delta$ of $\gg$.  Recall the central 
   element $z_{2\mu}$ introduced before Lemma 2.2.

     \begin{lemma}  Let $\gg,\gg^{\theta}$ be a symmetric pair of 
     type $DI$, case (i).   Assume further that 
     $\pi_{\Theta}=\{\alpha_{n-m+1},\alpha_{n-m+2}, \ldots, \alpha_{n-1},\alpha_{n}\}$
     where $n-1\geq m\geq 3$.
     Then  $${\cal P}_{B}(z_{2\omega_1})=a^{-1}(\sum_{\mu\in 
     W_{\Theta}\tilde\omega_1}q^{(\tilde\rho,2\mu)}\tau(2\mu)
     +\sum_{i=0}^{m-1}(q^{2i}+q^{-2i})$$ where 
     $$a=(\sum_{\mu\in 
     W_{\Theta}\tilde\omega_1}q^{(\tilde\rho,2\mu)}
     +\sum_{i=0}^{m-1}(q^{2i}+q^{-2i})).$$
     \end{lemma}
     
     \noindent
     {\bf Proof:}    First we compute ${\cal P}(z_{2\omega_1})$.   
     Note that $\omega_1$ is a minuscule weight.   In particular, 
     there does not exist $\beta\in P^+(\pi)$ with 
     $\beta<\omega_1$.   It follows from [L4, Lemma 5.1] that 
     $${\cal P}({z_{2\omega_1}})=a^{-1}(\sum_{\mu\in 
     W\omega_1}q^{(\rho,2\mu)}\tau(2\mu))$$ where $a=\sum_{\mu\in 
     W\omega_1}q^{(\rho,2\mu)}$.
     
     By [H, Section 13.2, Table 1],
     $$\omega_1=\alpha_1+\alpha_2+\dots+
     \alpha_{n-2}+(1/2)(\alpha_{n-1}+\alpha_{n}).$$   Recall  that we 
     are assuming 
    $\pi_{\Theta}=\{\alpha_{n-m+1},\alpha_{n-m+2},\dots
    \alpha_{n-1},\alpha_{n}\}$. It is straightforward to check using 
    [A] or 
    the list in
    [L2, Section 7] that $$\mu_{n-m}=\tilde\alpha_{n-m}= \alpha_{n-m}+
    \alpha_{n-m+1}+\cdots+\alpha_{n-2}+(\alpha_{n-1}+\alpha_{n})/2{\rm \ 
     \ }\eqno{(4.9)}$$
    $$\mu_i=\tilde\alpha_i=\alpha_i{\rm \ for \ }  1\leq 
    i\leq n-m-1\eqno{(4.10)}$$
    and $\tilde\alpha_j=0$ for $\alpha_j\in \pi_{\Theta}.$    Furthermore 
    the rank $t$ of $\Sigma$ equals $n-m$ and $\Sigma$ is of type 
    $B_{t}$.  It follows from [L4, Lemma 3.1] and [H, Section 13.2] that 
    $$ \tilde\omega_1=\eta_1=\sum_{1\leq i\leq 
    n-m}\tilde\alpha_i.\eqno{(4.11)}$$     From the description of the 
    fundamental weights corresponding to a root system of type 
    $B_{t}$, we see that the only restricted root in $P^+(\Sigma)$ 
    which is strictly less than $\eta_1$ is $0$.   Thus 
    the $W$ orbit of 
     $\omega_1$ can be written as a union of two sets
     $$S_1=\{\pm(\alpha_i+\cdots 
     +\alpha_{n-2}+(\alpha_{n-1}+\alpha_{n})/2)|
     1\leq i\leq n-m\}$$ and
     $$S_2= 
     \{\pm(\alpha_{i}+\alpha_{i+1}+\dots+\alpha_{n}-(\alpha_{n-1}+\alpha_{n})/2)|
     \ n-m+1\leq i\leq n-1\}$$
    It is straightforward to check  that $\mu=\tilde\mu$ and hence 
     $(\rho,\mu)=(\tilde\rho,\mu)$ for all $\mu\in 
     S_1$.  Moreover,  $S_1$ 
     corresponds to the $W_{\Theta}$ orbit of $\eta_1$.
     On the other hand, $S_2$ is a subset of $Q(\pi_{\Theta})$. 
     Hence, the 
     image under $\tilde{\ }$ of elements in the second set $S_2$ are all equal 
     to $0$.    
     
     We have $${\cal P}_B(z_{2\omega_1})=a^{-1}(\sum_{\mu\in 
     S_1}q^{(\tilde\rho,2\mu)}\tau(2\mu)+\sum_{\mu\in 
     S_2}q^{(\rho,2\mu)}).$$   The lemma now follows from the fact 
     that $$\sum_{\mu\in 
     S_2}q^{(\rho,2\mu)}=\sum_{i=0}^{m-1}(q^{2i}+q^{-2i}).\Box$$
     
     \medskip
     Continue  the assumption that $\gg,\gg^{\theta}$ 
     is a symmetric pair of type $DI$ case (i) and that 
     $\pi_{\Theta}=\{\alpha_{n-m+1},\alpha_{m-m+2},\dots, 
     \alpha_{n}\}$ with $n-1\geq m\geq 3$.   It follows that 
     $\omega_n$ is a fundamental minuscule weight in $P(\pi)$.  Furthermore,  
     $\tilde\omega_n=\eta_{t}$ ([L4, Lemma 3.2 and (3.5)]).  Since 
     $\Sigma$ is of type $B_t$, we see that $\eta_t$  is a
     minuscule fundamental weight with respect to the root system $\Sigma$.
     Let $\varphi_{ 2\eta_t}$  denote the  element of a $W_{\Theta}$ 
     invariant zonal spherical function associated to $2\eta_t$. Multiplying by a nonzero 
     scalar if necessary, we may assume that
     $\varphi_{2\eta_t}(\tau(\tilde\rho))=1$.   Since 
     $\eta_t$ is minuscule and $\varphi_{2\eta_t}$ is a $W_{\Theta}$ 
     invariant element of ${\cal C}[P(2\Sigma)]$, it follows that 
     $$\varphi_{2\eta_t}={{\sum_{\mu\in 
     W_{\Theta}\eta_t}z^{2\mu}}\over {\sum_{\mu\in 
     W_{\Theta}\eta_t}q^{(\tilde\rho,2\mu)}}}.\eqno{(4.12)}$$  
      Set 
     $c_{2\omega_1}={\cal L}(\tau(2\omega_1))$ and recall that $c_{2\omega_1}$ is 
     the unique element in $\check U^B$ which is also contained in 
     $\tau(2\omega_1)+(\adr B_+)\tau(2\omega_1)$.
     In the next lemma, we use 
     this zonal spherical function to distinguish between the $B$ 
     invariant element $c_{2\omega_1}$ and the central element 
     $z_{2\omega_1}$.
     
     \begin{lemma}  Suppose that $\gg,\gg^{\theta}$ is of type $DI$ 
     case (i).  Assume further that $\pi_{\Theta}=\{\alpha_{n-m+1},\alpha_{m-m+2},\dots, 
     \alpha_{n}\}$ with $n-1\geq m\geq 3$.   Then  ${\cal 
     P}_{B}(c_{2\omega_1})\neq {\cal P}_B(z_{2\omega_1})$.
     \end{lemma}
	 
     \noindent
     {\bf Proof:} By (1.2) we have 
     $$\varphi_{2\eta_t}( \tau(2\omega_1) 
     \tau(\tilde\rho))=z^{2\eta_t}({\cal P}_B(c_{2\omega_1})).$$
     Recall that $\Sigma$ is a root system of type $B_{t}$ with 
     $t=n-m$ and 
     set of positive simple roots 
     given by (4.9) and (4.10). By  [H, Section 13.2, Table 1], we 
    have
    $$\tilde\omega_n=\eta_{t}=1/2(\mu_1+2\mu_2+\cdots +
    (n-m)\mu_{{n-m}}).\eqno{(4.13)}$$  Furthermore, one checks that 
    $\sum_{\mu\in W_{\Theta}\eta_{n-m}}z^{2\mu}$ is equal to the 
    product $$\prod_{m\leq s\leq 
    n-1} 
    (z^{(\mu_{n-s}+\mu_{n-s+1}+\dots+\mu_{{n-m}})}-
     z^{-(\mu_{n-s}+
    \mu_{n-s+1}+\dots+ \mu_{{n-m}})}).\eqno{(4.14)}$$  
    Indeed, 
    writing (4.14) as a sum of elements of the form 
    $z^{\beta}$, we see that (4.14) is an element of the 
    set 
    $$z^{2\eta_{{t}}}+\sum_{\gamma\in Q^+(\Sigma)}{\bf 
    N}z^{2\eta_{{t}}-2\gamma}.$$  On the other hand,  
   it is  not hard to see that (4.14) 
    is $W_{\Theta}$ invariant.   Since  $\eta_{{t}}$ is minuscule 
    fundamental weight, 
    the desired   equality   follows. 
    
    By (4.9) and (4.10) we have
     $(\tilde\rho,\mu_{{t}})=(\rho,\mu_{{t}})=
     m$ while $(\tilde\rho,\mu_i)=(\rho,\alpha_i)=1$ for 
     $1\leq i\leq {n-m}-1$.  It follows from (4.12) and the    previous paragraph that 
     $$\varphi_{2\eta_{t}}=a^{-1}\prod_{m\leq s\leq 
    n-1} 
    (z^{(\mu_{n-s}+\mu_{n-s+1}+\dots+\mu_{{n-m}})}-
     z^{-(\mu_{n-s}+
    \mu_{n-s+1}+\dots+ \mu_{{n-m}})})$$ where
  $$a=\prod_{m\leq s\leq 
    n-1} 
    (q^{s}+q^{-s}).$$
     Now  (4.9) and (4.10) and the fact that $\pi$ generates a root 
     system of type $D_n$ ensure that $( \omega_1,\mu_i)=
     \delta_{i1}(\omega_1,\alpha_1)=\delta_{ij}.$  Hence 
     $$\eqalign{\varphi_{2\eta_t}(\tau(2\omega_1+\tilde\rho))&=
  a^{-1} (q^{n+1}+q^{-n-1}) \prod_{m\leq s\leq 
    n-2} 
    (q^{s}+q^{-s})\cr&=
    (q^{n+1}+q^{-n-1})(q^{n-1}-q^{-n+1})^{-1}.\cr}$$

      Furthermore, it is straightforward to 
    check using (4.11) that the $W_{\Theta}$ orbit of $\eta_1$
    is the set $\{\pm(\mu_{n-s}+\mu_{n-s+1}+\cdots +\mu_{{n-m}}|\ 
    m\leq s\leq n-1\}$.  By 
    Lemma 4.7, $b{\cal P}_B(z_{2\omega_1})-\sum_{i=0}^{m-1}(q^m+q^{-m})$ equals
     $$ \sum_{m\leq s\leq n-1} q^{2s}\tau(
     2(\mu_{n-s}+\cdots +\mu_{n-m}))+
     q^{-2s}\tau(-2(\mu_i+\cdots +\mu_{n-m}))
     $$  where $b=(\sum_{m\leq s\leq 
     n-1}q^{2s}+q^{-2s})+\sum_{i=0}^{m-1}(q^{2i}+q^{-2i})$. 
     Hence using (4.13) we have
     $$z^{2\eta_t}({\cal P}_B(z_{2\omega_1}))=
   b^{-1}  (q^{2n}+q^{-2n}-(q^{2(n-1)}+q^{-2(n-1)})+b) .$$
    Thus ${\cal P}_B(z_{2\omega_1})={\cal P}_B(c_{2\omega_1})$ implies that 
    $$(q^{n+1}+q^{-n-1})b
    =(q^{2n}+q^{-2n}-(q^{2(n-1)}+q^{-2(n-1)})+b)(q^{n-1}+q^{-n+1}).\eqno{(4.15)}$$
  Note that the coefficient of $q^{n-1}$ is $0$ in the left hand side of 
  (4.15).   On the other hand, the coefficient of $q^{n-1}$ in the 
  right hand side of (4.15) is $2$.  This 
   contradiction proves the lemma.  $\Box$
   
   \medskip
   
   For Cases (iii) and (iv) below, it is convenient to use different 
   notation for the fundamental restricted weights  
    in $P^+(\Sigma)$.  
   Given $\alpha_i\notin\pi_{\Theta}$, let $\omega_i'$ denote the 
   fundamental weight associated to the simple restricted root 
   $\tilde\alpha_i$.    
   
   \bigskip
   \noindent
   {\bf Case (iii):}
   We now consider the three exceptional types of symmetric pairs, 
   EIV, EVII, and EIX.    In each case, we have that the 
   set of fundamental weights in $P^+(\Sigma)$ not contained in the image 
   of $P^+(\pi)$
     under $\tilde{\ }$ is precisely the set
   $\{ \omega_1', \omega_6'\}$.  Furthermore, one checks that 
  $\omega_6'$ is the unique nonzero element of $P^+(\Sigma)$ strictly 
  less than $\tilde\omega_1=2\omega_1'$ while $\omega_1'$ is the unique nonzero 
  element of $P^+(\Sigma)$ strictly less than $\tilde\omega_n$.  
  Here $\tilde\omega_n=\omega_n'$ in the latter two cases (i.e. $n=7$ 
  or $n=8$) while $\tilde\omega_n=2\omega_6'$ when 
  $\gg,\gg^{\theta}$ is of type EIV.
  
  Theorem 4.1 follows from Lemma 4.3 and  the next lemma for Case (iii) symmetric pairs.  
  \begin{lemma}
      Assume that $\gg,\gg^{\theta}$ is of type EIV, EVII, or 
     EIX.   Then there exists $f\in (\adr U)\tau(2\omega_1)$ 
      such that $$f\in\tau(2\omega_6')+
      (B\check T_{\Theta})_+\check U$$
      and $g\in (\adr U)\tau(2\omega_n)$ such that
      $$g\in\tau(2\omega_1')+
      (B\check T_{\Theta})_+\check U.$$
 \end{lemma}
 
 \noindent
 {\bf Proof:}  Let $\pi'$ be the subset of $\pi$ equal to 
 $\{\alpha_i|\ i>1\}$.  Note that $\pi'$ generates a root system of 
 type $D_{n-1}$.   Furthermore,   
  $\alpha_{n-i}$ is the $(i+1)^{th}$ simple root in this root system with respect to 
 the ordering of the simple roots given in [H, 
 Chapter III].  Now $\Theta$ 
 restricts to an involution on the root system generated by $\pi'$.
 Moreover, $\pi'\cap \pi_{\Theta} 
 =\pi_{\Theta}=\{\alpha_2,\alpha_3,\alpha_4,\alpha_5\}$.   Let $\gr$
 be the semisimple Lie subalgebra 
 of $\gg$ generated by the positive and negative root vectors 
 corresponding to the simple roots in $\pi'$.   The Lie algebra $\gr$ has 
 rank $n-1$ and $\gr, 
\gr^{\theta}$ is a symmetric pair of type DI, case (i).      
 We write $U_q(\gr)$ for the  quantized enveloping algebra of $\gr$ 
 indentified in the obvious way with a subalgebra of $U$. Let $\check 
 U_q(\gr)$ denote the simply connected quantized enveloping algebra 
 of $\gr$. Let $\Sigma'$ be the restricted root system associated to the 
 symmetric pair $\gr,\gr^{\theta}$ and let $W'_{\Theta}$ denote the 
 corresponding restricted Weyl group.
 
Let $\nu_n$ denote the fundamental weight associated to the simple 
root $\alpha_n$ considered as an element in the root system 
associated to $\pi'$. Let $z'_{2\nu_{n}}$ be the unique central element 
of   $\check  U_q(\gr)$ 
such that $$z'_{2\nu_n}\in\tau(2\nu_{n})+(\adr 
U_q(\gr)_+)\tau(2\nu_n).$$  
Similarly, Let $c'_{2\nu_{n}}$ be the unique ${(U_q(\gr)\cap B)}$ 
invariant  element 
of   $ \check U_q(\gr)$ 
such that  $$c'_{2\nu_n}\in\tau(2\nu_{n})+(\adr (U_q(\gr)\cap 
B)_+)\tau(2\nu_n).$$
  Note that both $z'_{2\nu_{n}}$ and $c'_{2\nu_n}$ are elements of 
$(\adr U_q(\gr))\tau(2\nu_n)$.  

Now $\Sigma'$ is a root system of type $B_{n-5}$. Moreover, 
if we order the roots of $\Sigma'$ as in [H], then 
$\tilde\alpha_n$ corresponds to the first simple root of $\Sigma'$. 
By (4.11), $\tilde\nu_n$ is the fundamental weight corresponding to 
$\tilde\alpha_n$. It follows that $\tilde\nu_n$ is a pseudominuscule 
weight in $P(\Sigma')$.  Hence, Lemma 2.2 implies that both    ${\cal 
P}_{B\cap 
U_q(\gr)}(z'_{2\nu_n})$ and 
${\cal P}_{B\cap 
U_q(\gr)}(c'_{2\eta_n})$ are linear combinations of $1$ and 
$\sum_{\gamma\in W'_{\Theta}\nu_n}q^{(\rho,2\gamma)}\tau(2\gamma).$  

 Lemma 4.8 ensures 
that ${\cal P}_{B\cap U_q(\gr)}(z'_{2\nu_n})\neq {\cal P}_{B\cap 
U_q(\gr)}(c'_{2\nu_n})$. 
Hence there is a linear combination $X$ of  $ z'_{2\nu_n}$ and 
$ c'_{2\nu_n}$ such that ${\cal P}_{B\cap 
U_q(\gr)}(X)=1$.   It follows that 
${\cal P}_{B\cap U_q(\gr)}(X-1)=0$.  Thus by [L4, Corollary 4.2],  $X-1\in (B\check T_{\Theta}\cap 
\check U_q(\gr))_+U$.   In particular, 
$$X\in (\adr U_q(\gr))\tau(2\nu_n){\rm \ and \ }X\in 1+ (B\check 
T_{\Theta})_+\check U.\eqno{(4.16)}$$
Since $\alpha_n$ is the first root in the root system of type 
$D_{n-2}$ generated by $\pi'$, we have  
$$\nu_n=1/2(\alpha_2+\alpha_3)+\alpha_4+\dots+\alpha_n.$$
Note that $(\omega_n-\nu_n,\alpha_i)=0 $  for all $\alpha_i\in \pi'$.
On the other hand $$(\omega_n-\nu_n,\alpha_1)=(-\nu_n,\alpha_1)=
-(\alpha_3,\alpha_1)/2=(\omega_1,\alpha_1)/2.$$ It follows that 
$\omega_n=\nu_n+\omega_1/2.$    Now $\tau(\omega_1)$ 
commutes with elements of $U_q(\gr)$.   Hence
$$(\adr U_q(\gr))\tau(2\omega_n)=((\adr 
U_q(\gr))\tau(2\nu_n))\tau(\omega_1).\eqno{(4.17)}$$
It follows from (4.16) and (4.17) that 
$$X\tau(\omega_1) \in (\adr U)\tau(2\omega_n){\rm \ and \ }X\tau(\omega_1)\in 
\tau(\tilde\omega_1)+ (B\check 
T_{\Theta})_+\check U.$$
The second assertion now follows from the fact that 
$\tilde\omega_1=2\omega'_1$ ([L4, Lemma 3.1]).   The first assertion is proved in 
exactly the same way where we replace $\pi'$ with the set 
$\{\alpha_i|\ i<5\}$. $\Box$

\bigskip
\noindent
{\bf Case (iv):}
We next turn our attention to the symmetric pair of type EIII.    
In this case, $\tilde\alpha_1=(\alpha_1+\alpha_3+\alpha_4+\alpha_5+\alpha_6)/2$ and 
$\tilde\alpha_2= \alpha_2+\alpha_4+(\alpha_3+\alpha_5)/2$.
Moreover,  both $\tilde\alpha_1$ and $2\tilde\alpha_1$ are elements 
of $\Sigma$.  In particular, the restricted root system $\Sigma$
is nonreduced of type   BC$_2$. Now $\Sigma$ contains a root system 
with set of positive simple roots $\{2\tilde\alpha_1,\tilde\alpha_2\}$ 
of type  B$_2$. Here $\tilde\alpha_2$ is the short simple root and 
$2\tilde\alpha_1$ is the long simple root. As explained in [L4, Section 3], 
the fundamental weight $\omega_1'$ associated to $\tilde\alpha_1$ 
satisfies 
$(\omega_1',2\tilde\alpha_1)=(2\tilde\alpha_1,2\tilde\alpha_1)/2$.  In 
particular, the weight lattice $P(\Sigma)$ is the same as the weight 
lattice of the underlying root system of type  B$_2$. 
Thus 
$\omega_1'=2\tilde\alpha_1+\tilde\alpha_2$ and 
$\omega_2'=\tilde\alpha_1+\tilde\alpha_2$. It is 
straightforward to check that
$\tilde\omega_1=2\tilde\alpha_1+\tilde\alpha_2$ and 
$\tilde\omega_2=2\tilde\alpha_1+2\tilde\alpha_2$.    Note that $\omega_1'\in 
\widetilde{P^+(\pi)}$ while $\omega_2'$ is not an element of this 
set. (This last fact is also used in the proof of [L4, Theorem 2.5]).  

By Lemma 4.3,  we have that 
$\tau(2\omega_1')+(\tau(2\omega_1)-\tau(2\omega_1'))\in 
F_r(\check U)$ and $\tau(2\omega_1)-\tau(2\omega_1')\in (B\check T_{\Theta})_+\check U$. Thus 
Theorem 4.1 for this last case follows from the next lemma.

\begin{lemma} Assume that $\gg,\gg^{\theta}$ is of type EIII. 
    Then there exists $f\in F_r(\check U)$ 
      such that $$f\in\tau(2\omega_2')+
      (B\check T_{\Theta})_+\check U.$$
    \end{lemma}
    
    \noindent
    {\bf Proof:}  Let $\pi'$ be the subset of $\pi$ equal to 
    $\{\alpha_2,\alpha_3,\alpha_4,\alpha_5\}$,   Note that $\pi'$ 
    generates a root system of type D$_4$.   Furthermore, $\alpha_i$ 
    is the $(i-1)^{th}$ simple root with respect to the ordering of 
    the simple roots given in [H].   Now $\Theta$ restricts to an 
    involution on the root system generated by $\pi'$.   Note that 
    $\pi'\cap \pi=\{\alpha_3,\alpha_4,\alpha_5\}$. Let $\nu_2$ denote 
    the fundamental weight associated to the simple root $\alpha_2$ 
    with respect to the root system generated by $\pi'$.
  We have $\nu_2=(\alpha_3+\alpha_5)/2+\alpha_4+\alpha_2$.   Now
$(\omega_2-\nu_2,\alpha_i)=(\omega_i,\alpha_i)/2$ for $i=1$ and 
$i=6$.  Also, $(\omega_2-\nu_2,\alpha_i)=0$ for $i\notin\{1,6\}$.  
Therefore $$\omega_2=\nu_2+(\omega_1+\omega_6)/2.$$
By [L4, Lemma 3.1], $\tilde\omega_1+\tilde\omega_6=2\omega_1'$.  The 
rest of the argument follows  as in the proof of Lemma 4.9. $\Box$

\section{Appendix: Commonly used  notation}

\noindent
Here is a list of notation defined in Section 1 (in the following 
order):

\noindent
${\bf C}$, ${\bf Q}$,
 ${\bf R}$,
$q$,
${\cal C}$, ${\cal R}$,
$Q(\Phi)$, $P(\Phi)$, $Q^+(\Phi)$, $P^+(\Phi)$, $\gg$, $\gn^-$, 
$\gh$, $\gn^+$, $\Delta$, $\pi=\{\alpha_1,\dots, 
\alpha_n\}$, $(\ , \ )$, $\leq $, $\theta$, $\gg^{\theta}$, $\Theta$,
$\pi_{\Theta}$, $\tilde\alpha$, $\Sigma$, $U$, $x_i$, $y_i$, $t_i^{\pm 
1}$, $T$, $U^0$, $U^+$, $G^-$, $\tau$, $\check U$, $\check U^0$, $A_+$, ${\cal B}$, ${\cal C}[G]$, 
${\cal C}[H]$, $\adr$, $\ad$, $\check U^B$, $\check{\cal A}$, $\check 
T_{\Theta}$, ${\cal M}$, $N^+$, ${\cal P}_B$, ${\cal A}$, $*$, 
 $\varphi_{\lambda}$,  ${\cal X}$, ${\cal C}(Q(\Sigma)){\cal A}$.

 \medskip
 \begin{tabbing}
	 \noindent
 \=Defined in Section 2:\=\\
 \>$Z(\check U)$\> the center of $\check U$\\
 \>$\rho$\> half sum of positive roots in $\Delta$\\
 \>$w.q^{(\rho,\lambda)}\tau(\lambda)$\>$q^{(\rho,w\lambda)}\tau(w\lambda)$\\
 \>$\hat m(2\eta)$\>$\sum_{\gamma\in 
 W_{\Theta}\eta}q^{(\rho,2\gamma)}\tau(2\gamma)$\\
 \>${\rm tip}$\> highest degree  term w.r.t. ``deg'' 
 degree function\\
 \>${\cal C}[[Q(\Sigma)]]$\>power series ring in the $z^{-\gamma}$ 
 for $\gamma\in Q^+(\Sigma)$\\
 \>${\cal N}_{\eta}$\> see (2.3)\\
 \>${\cal N}_{\eta}^+$\> see (2.3)\\
 \>$w_o'$\> the longest element of $W_{\Theta}$\\
 \>$z_{2\mu}$\> unique central element in $\tau(2\mu)+(\adr 
 U_+)\tau(2\mu)$\\
 \>${\rm top}$\> highest degree  term w.r.t.
 ``odeg'' degree function\\
 \>\>\\
 \=Defined in Section 3:\=\\
 \>$\phi$\> Hopf algebra automorphism of $U$ (see (3.1))\\
 \>$T_{\Theta}$\>$\{\tau(\beta)|\Theta(\beta)=\beta$ and $\beta\in 
 Q(\pi)\}$\\
 \>$y_it_i+d_i\tilde\theta(y_i)t_i+s_it_i$\> one of the  generators 
 of $B$\\
 \>$\tilde\theta$\> lift of the involution $\theta$ to $U$\\
 \>$\sigma$\> the antipode of $U$\\
 \>${\it\Delta}$\> the coproduct of $U$\\
 \>$\chi$\> see definition preceding Lemma 3.2\\
 \>$N^-$\> algebra generated by $(\ad {\cal M}\cap G^-)[y_it_i,
  \alpha_i\notin\pi_{\Theta}]$\\
 \>${\cal P}'_B$\> projection onto ${\cal C}[\check {\cal A}]$ defined 
 using (3.2)\\
 \>$\kappa$\> antiautomorphism of $U$ preserving $B$\\
 \>$L(\lambda)$\> highest weight simple $U$ module of highest weight 
  $\lambda$\\
 \>$F_r(\check U)$\> locally finite part of $\check U$ w.r.t $\adr$\\
 \>${\cal L}$\> see Definition 3.4\\
 \>$g_{\lambda}$\> zonal spherical function at $\lambda$ associated to 
 $\chi^{-1}(B),B$\\
 \>$\xi_{\lambda}$\> nonzero $B$ invariant vector of $L(\lambda)$\\
 \>$\zeta_{\lambda}^*$\> nonzero $\chi^{-1}(B)$ invariant vector of 
 $L(\lambda)^*$\\
 \>$\hat N^+$\>$\sum_{\gamma\in Q^+(\pi)}N^+_{\gamma}\tau(-\gamma)$\\
 \>$\hat N^-$\>$\sum_{\gamma\in Q^+(\pi)}N^-_{-\gamma}\tau(-\gamma)$\\
 \>$T'_{\leq}$\> $\{\tau(-\gamma)|\gamma\in Q^+(\pi)$ and 
 $\tilde\gamma\in P(2\Sigma)\}$\\
 \>$S_{\beta,r}$\> sum of weight spaces $S_{\beta'}$ with 
 $\tilde\beta'=\tilde\beta$\\
 \>$G^+$\> algebra generated by $x_it_i^{-1}$, $1\leq i\leq n$\\
 \>$U^-$\> algebra generated by $y_i$, $1\leq i\leq n$\\
 \>${\cal A}_{\leq}$\> $\{\tau(-\gamma)|\ \gamma\in 
 Q^+(\Sigma)\cap P(2\Sigma)\}$\\
 \>\>\\
  \=Defined in Section 4:\=\\
 \>$t$\> rank of the restricted root system $\Sigma$\\
 \>$\mu_1,\dots, \mu_t$\> simple roots in $\Sigma$\\
 \>$\eta_1,\dots, \eta_t$\> fundamental weights in $P^+(\Sigma)$\\
 \>$\pi_i$\> $\{\alpha_j|\ 1\leq j\leq 2i+1\}$\\
 \>$\gg_i$\> semisimple Lie subalgebra of $\gg$ with simple roots 
 $\pi_i$\\
 \>$U_i$\> subalgebra of $U$ equal to $U_q(\gg_i)$\\
 \>$U_0$\>${\cal C}$\\
 \>$U_t$\>$U$\\
 \>$\check U_i$\>simply connected quantized enveloping algebra of 
 $\gg_i$\\
 \>$\check U^i$\> algebra generated by $U_i$ and $\check T$\\
 \>$\omega_j$\> fundamental weight corresponding to $\alpha_i$ in 
 $P^+(\pi)$\\
 \>${\cal P}$\> ordinary Harish-Chandra projection of $\check U$ 
 onto $\check U^0$\\
 \>$W$\> Weyl group of root system $\Delta$\\
 \>$c_{2\omega_1}$\>${\cal L}(\tau(2\omega_1))$\\
 \>$\omega_i'$\> restricted fundamental weight associated to $\tilde\alpha_i$\\

 \end{tabbing}

 \centerline{REFERENCES}

\bigskip
\bigskip
\noindent
[A] S. Araki, On root systems and an infinitesimal 
classification of irreducible symmetric spaces, {\it Journal of 
Mathematics, Osaka City University} {\bf 13} (1962), no. 1, 1-34.

\medskip
\noindent
[B] P. Baumann, On the center of quantized enveloping algebras,
{\it Journal of Algebra}  {\bf 203} (1998), 244-260.

\medskip
\noindent [DN] M.S. Dijkhuizen and M. Noumi, A family of quantum
projective spaces and related $q$-hypergeometric orthogonal
polynomials, {\it Transactions of the A.M.S.} {\bf 350} (1998), no.
8, 3269-3296.

\medskip
\noindent [DS] M.S. Dikhhuizen and J.V. Stokman, Some limit transitions between BC type orthogonal polynomials interpreted on
quantum complex Grassmannians, {\it Publ. Res. Inst. Math. Sci. 35} (1999), 451-500.

\medskip\noindent [H]   J.E. Humphreys, {\it Introduction to Lie Algebras
and Representation Theory}, Springer-Verlag, New York
(1972). 

\medskip\noindent[HC] Harish-Chandra, Spherical functions on a 
semisimple Lie group, I, {\it American Journal of Mathematics}  {\bf 80} 
(1958) 241-310.

\medskip \noindent [JL] A. Joseph and G. Letzter, Local
finiteness of the adjoint action for quantized enveloping
algebras, {\it Journal of Algebra} {\bf 153} (1992), 289 -318.

\medskip\noindent  [Jo] A. Joseph, {\it Quantum Groups and Their Primitive
Ideals}, Springer-Verlag, New York (1995).

\medskip\noindent [KS] E. Koelink, J. Stokman,   Fourier transforms 
on the quantum SU(1,1) group.  With an appendix by Mizan Rahman {\it Publ. 
Res. Inst. Math. Sci.} {\bf 37} (2001), no. 4, 621-715. 

\medskip\noindent [K] A.A. Kirillov, Jr., Lectures on affine Hecke 
algebras and Macdonald's conjectures, {\it Bulletin of the American 
Mathematical Society} {\bf 34} (1997), No. 3, 251-292.

\medskip \noindent [L1] G. Letzter, Coideal subalgebras and quantum 
symmetric pairs, In: {\it New Directions in Hopf Algebras, MSRI publications}
{\bf 43}, Cambridge University Press (2002), 117-166.

\medskip\noindent [L2]  G. Letzter, Quantum symmetric pairs and their 
zonal spherical functions, {\it Transformation Groups} {\bf 8} 
(2003), no. 3. 261-292.

\medskip\noindent [L3]  G. Letzter, Quantum  zonal spherical  
functions and Macdonald polynomials, Advances in Mathematics, in
press (corrected proof available online).
 
\medskip\noindent [L4] G. Letzter, Invariant differential operators 
for quantum symmetric spaces, I, preprint (arXiv:math.QA/0406193).

\medskip
\noindent
[M] I.G. Macdonald,
Orthogonal polynomials associated with root systems, {\it S\'eminaire Lotharingien de Combinatoire} 
{\bf 45} (2000/01) 40 pp.

\medskip
\noindent [N] M. Noumi, Macdonald's symmetric polynomials as zonal
spherical functions on some quantum homogeneous spaces, {\it
Advances in Mathematics} {\bf 123} (1996), no. 1, 16-77.

\medskip
\noindent [NDS] M. Noumi, M.S. Dijkhuizen, and T. Sugitani, Multivariable Askey-Wilson polynomials and quantum complex
Grassmannians, {\it Fields Institute Communications} {\bf 14} (1997), 167-177.

\medskip
\noindent [NS] M. Noumi and T. Sugitani, Quantum symmetric spaces
and related q-orthogonal polynomials, in: {\it Group Theoretical
Methods in Physics (ICGTMP)} (Toyonaka, Japan, 1994), World
Science Publishing, River Edge, New Jersey (1995), 28-40.

\medskip
\noindent [S] T. Sugitani, Zonal spherical functions on quantum Grassmann manifolds, {\it J. Math. Sci. Univ. Tokyo}
{\bf 6} (1999), no. 2, 335-369.

\end{document}